\documentclass[12pt,letterpaper]{amsart}

\usepackage{graphicx,amscd}
\newtheorem{theorem}{Theorem}[section]
\newtheorem{proposition}[theorem]{Proposition}
\newtheorem{lemma}[theorem]{Lemma}
\newtheorem{corollary}[theorem]{Corollary}
\newtheorem*{claim}{Claim}

\theoremstyle{remark}

\begin{document}

\title{Reducible and toroidal Dehn fillings with distance 3}

\author[S. Kang]{Sungmo Kang}
\address{Department of Mathematics, University of Texas at Austin, 1 University Station C1200 Austin, TX
78712-0257} \email{skang@math.utexas.edu}

\maketitle

\begin{abstract}
If a simple $3$-manifold $M$ admits a reducible and a toroidal Dehn
filling, the distance between the filling slopes is known to be
bounded by three. In this paper, we classify all manifolds which
admit a reducible Dehn filling and a toroidal Dehn filling with
distance 3.
\end{abstract}

{\em Keywords}: Dehn filling; Reducible; Toroidal; Tangle filling.\

\maketitle

\section{Introduction}
Let $M$ be a compact connected orientable $3$-manifold with a torus
boundary component $\partial_0 M$ and $r$ a {\em slope}, the isotopy
class of an essential simple closed curve, on $\partial_0 M$. The
manifold obtained by {\em $r$-Dehn filling} is defined to be
$M(r)=M\cup V$, where $V$ is a solid torus glued to $M$ along
$\partial_0 M$ so that $r$ bounds a disk in $V$.

We say that $M$ is {\em simple} if it contains no essential sphere,
torus, disk or annulus. For a pair of slopes $r_1$ and $r_2$ on
$\partial_0 M$, the distance $\Delta (r_1,r_2)$ denotes their
minimal geometric intersection number. For a simple manifold $M$, if
both $M(r_1)$ and $M(r_2)$ fail to be simple, then the upper bounds
for $\Delta(r_1,r_2)$ have been established in various cases. In
particular, for the case of reducible and toroidal Dehn fillings, Oh
\cite{O1} and independently Wu \cite{W3} showed that for a simple
manifold $M$, if $M(r_1)$ is reducible and $M(r_2)$ is toroidal then
$\Delta(r_1,r_2)\le 3$. Eudave-Mu\~{n}oz and Wu \cite{EW} have
explicitly described an infinite family of manifolds $M$ realizing
$\Delta(r_1, r_2) = 3$ in terms of tangle arguments. More
specifically, $M$ is the double branched cover of one of the
tangles illustrated in Figure 1 where $p$ is an integer no less than
3. However, $p$ can be extended to any integer no less than 2.

The goal of this paper is to show that the manifolds $M$ described
above are the only examples of simple manifolds admitting a
reducible Dehn filling and a toroidal Dehn filling with distance 3.
The main results of this paper are the following.
\begin{figure}[tbh]
\includegraphics{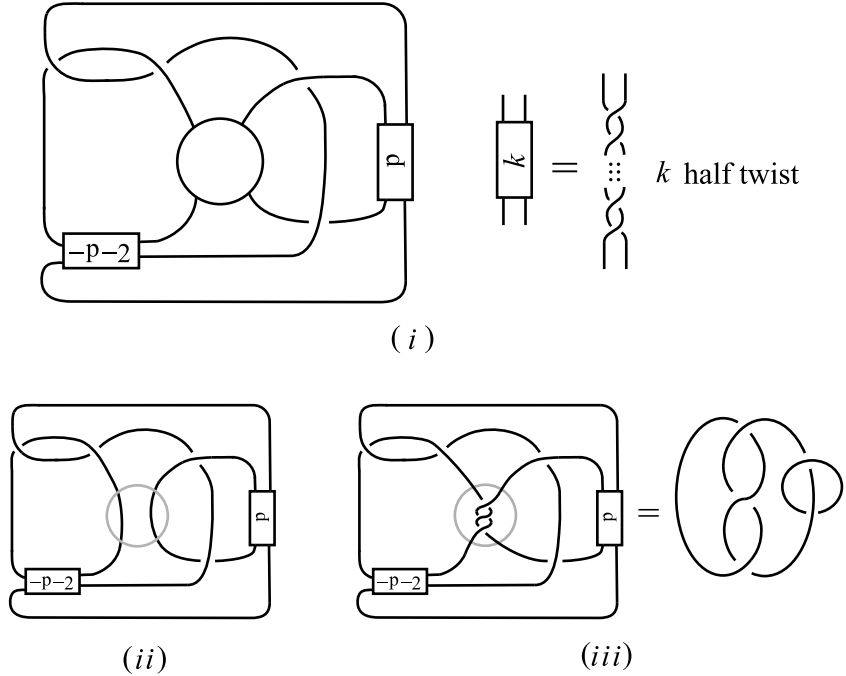}\caption{The tangles and the $1/0$- and $1/3$-rational tangle fillings used in Theorem 1.1.}\label{f:jj}
\end{figure}

\begin{theorem} \label{main1}
Let $M$ be a simple $3$-manifold with a torus boundary component
$\partial_0 M$ such that $M(r_1)$ is reducible, $M(r_2)$ is toroidal
and $\Delta(r_1,r_2)=3$. Then M is the double branched cover of one
of the tangles with $p \geq 2$ as described in Figure $1$ $(i)$.
$M(r_1)$ corresponding to the $1/3$-rational tangle filling is a
connected sum of two lens spaces $L(3, 1)$ and $L(2, 1)$, and
$M(r_2)$ corresponding to the $1/0$-rational tangle filling is
toroidal and is not Seifert fibered. See Figure $1$ $(ii), (iii)$.
\end{theorem}

In order to prove Theorem \ref{main1}, we need the following.

\begin{theorem} \label{main}
Let $M$ be a simple $3$-manifold as in Theorem $\ref{main1}$. Then
$M(r_2)$ is the union of $M_1$ and $M_2$, where $M_i$, $i=1, 2$, is
a Seifert fibered space over a disk with two exceptional fibers,
along an incompressible torus. Furthermore, the fibers of $M_1$ and
$M_2$ intersect exactly once on the incompressible torus i.e.
$M(r_2)$ is not Seifert fibered.
\end{theorem}

In \cite{Lee}, Lee showed that if $M(r_1)$ is reducible and $M(r_2)$
contains a Klein bottle, then $M$ is uniquely determined. That is,
$M$ is homeomorphic to $W(6)$, which is obtained from the exterior
$W$ of the Whitehead link by performing Dehn filling on one of its
boundary components with slope 6 under the standard
meridian-longitude coordinates. Let $M=W(6)$. Then it is well-known
that $M(1)$ is reducible, and $M(4)$ is toroidal and contains a
Klein bottle. Moreover, using the fact that $W$ is strongly
invertible, we can easily see that $M$ is the double branched cover
of the tangle illustrated in Figure 1 ($i$) with $p=2$, say $T$, and
$M(1)$ and $M(4)$ are the double branched covers of $T$ with the
1/3- and 1/0-rational tangles filled respectively. Notice that the
toroidal Dehn filling $M(4)$ is the union of $M_1$ and $M_2$, where
$M_i$, $i=1, 2$, is a Seifert fibered space over a disk with two
exceptional fibers, along an incompressible torus and the fibers of
$M_1$ and $M_2$ intersect exactly once on the incompressible torus
i.e. $M(4)$ is not Seifert fibered. We have shown that if $M(r_2)$
contains a Klein bottle, then Theorem \ref{main1} and Theorem
\ref{main} hold.

Therefore, throughout this paper we may assume that $M(r_2)$ does not
contain a Klein bottle. Then Theorem \ref{main1} and Theorem
\ref{main} follow immediately from  Theorem
\ref{main2} and Theorem \ref{23} respectively.

As an immediate consequence of Theorem \ref{main}, we have the
following.
\begin{corollary} \label{cor}
Let $M$ be a simple manifold with a torus boundary component. If
$M(r_1)$ is reducible and $M(r_2)$ is a toroidal Seifert fibered
manifold, then $\Delta(r_1,r_2)\le 2$.
\end{corollary}

It is unknown whether or not the upper bound 2 in Corollary
\ref{cor} is the best possible.

\section{The intersection graphs}

From now on we assume that $M$ is a simple $3$-manifold with a torus
boundary component $\partial_0 M$ and that $r_1$ and $r_2$ are
slopes on $\partial_0 M$ such that $M(r_1)$ is reducible and
$M(r_2)$ contains an essential torus and $\Delta (r_1,r_2)=3$.

Over all reducing spheres in $M(r_1)$ which intersect the attached
solid torus $V_1$ in a family of meridian disks, we choose a
$2$-sphere $\widehat{F}_1$ so that $F_1=\widehat{F}_1 \cap M$ has the
minimal number, say $n_1$, of boundary components. Similarly let
$\widehat{F}_2$ be an essential torus in $M(r_2)$ which intersects
the attached solid torus $V_2$ in a family of meridian disks, the
number of which, say $n_2$, is minimal over all such surfaces and
let $F_2=\widehat{F}_2\cap M$. Let $u_1,u_2,\ldots, u_{n_1}$ be the
disks of $\widehat{F}_1\cap V_1$, labeled as they appear along
$V_1$. Similarly let $v_1,v_2,\ldots,v_{n_2}$ be the disks of
$\widehat{F}_2\cap V_2$. Then $F_1$ is an essential planar surface,
and $F_2$ is an essential punctured torus in $M$. We may assume that
$F_1$ and $F_2$ intersect transversely and the number of components
of $F_1\cap F_2$ is minimal over all such surfaces. Then no circle
component of $F_1\cap F_2$ bounds a disk in either $F_1$ or $F_2$
and no arc component is boundary-parallel in either $F_1$ or $F_2$.
In this section, the subscripts $i$ and $j$ will denote 1 or 2 with
the convention that if they both appear, then $\{i,j\}=\{1,2\}$.
The components of $\partial F_i$ are numbered $1,2,\ldots ,n_i$
according to the labels of the corresponding disks of
$\widehat{F}_i\cap V_i$. We obtain a graph $G_i$ in $\widehat{F}_i$
by taking as the (fat) vertices of $G_i$ the disks in
$\widehat{F}_i\cap V_i$ and as the edges of $G_i$ the arc components
of $F_1\cap F_2$ in $F_i$. Each endpoint of an edge of $G_i$ has a
label, that is, the label of the corresponding component of
$\partial F_j$, $i\ne j$. Since each component of $\partial F_i$
intersects each component of $\partial F_j$ in $\Delta$ points, the
labels $1,2,\ldots ,n_j$ appear in order around each vertex of $G_i$
repeatedly  $\Delta$ times.

For a graph $G$, the {\em reduced graph} $\overline{G}$ of $G$ is
defined to be the graph obtained from $G$ by amalgamating each
family of mutually parallel edges into a single edge. For an edge
$\alpha$ of $\overline{G}$, the {\em weight} of $\alpha$, denoted by
$w(\alpha)$, is the number of edges of $G$ in $\alpha$.

 Orient all components of $\partial F_i$ so that they are
mutually homologous on $\partial_0 M$, $i=1,2$. Let $e$ be an edge
in $G_i$. Since $e$ is a properly embedded arc in $F_i$, it has a
disk neighborhood $D$ in $F_i$ with $\partial D=a\cup b\cup c \cup
d$, where $a$ and $c$ are arcs in $\partial F_i$ with induced
orientation from $\partial F_i$. On $D$, if $\partial D$ can be
oriented so that $a$ and $c$ have the same orientation as that
induced from $\partial F_i$, then $e$ is called {\em positive},
otherwise {\em negative}. See Figure 2. Then we have the parity
rule.

{\em An edge is positive on one graph if and only if it is negative
on the other graph.}

\begin{figure}[tbh]
\includegraphics{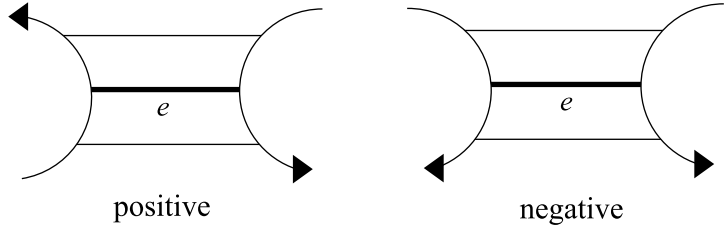}\caption{Orientations on an edge.}\label{f:jj}
\end{figure}

Orient the core of $V_i$. Then we can give a sign to each vertex of
$G_i$ according to the sign of its intersection with the core of
$V_i$.  Two vertices (possibly equal) of $G_i$ are called {\em
parallel} if they have the same sign, otherwise {\em antiparallel}.
A positive edge connects parallel vertices while a negative one
connects antiparallel vertices. Let $G^+_i$ denote the subgraph of
$G_i$ consisting of all the vertices and all the positive edges of
$G_i$.

Let $G$ be $G_1$ or $G_2$ and let $x$ be a label of $G$. An {\em
$x$-edge} is an edge of $G$ with the label $x$ at one endpoint. An {\em
$x$-cycle} is a cycle of positive $x$-edges which can be oriented so
that the tail of each edge has the label $x$. A cycle in $G$ is a {\em
Scharlemann cycle} if it bounds a disk face, and the edges in the
cycle are all positive and have the same label pair. If the label
pair is $\{x,y\}$, then we refer to such a Scharlemann cycle as an
{\em $(x,y)$-Scharlemann cycle}. In particular, a Scharlemann cycle
of length $2$ is called an {\em $S$-cycle}.

\begin{lemma} \label{properties}
If $G_i$ has a Scharlemann cycle, then $\widehat{F}_j$ is separating.
\end{lemma}
\begin{proof}
This follows from \cite[Lemma 1.2]{W3}.
\end{proof}

\section {Binary faces and the subgraph $\Lambda$ of $G_1^+$}
In this section, we will define binary faces in $G_1^+$ and see some
properties of binary faces. Since the arguments here are exactly the
same as in \cite{GL4} and in \cite{Lee} at the same time, we
use the terminologies and several Lemmas without proof described in
\cite{GL4} and in \cite{Lee}. We may assume that $M(r_2)$ is
irreducible by \cite{GLu3}.

\begin{lemma} \label{anti}
$G_2$ has exactly two vertices and the two vertices of $G_2$ are antiparallel.
\end{lemma}
\begin{proof}
By Theorem 1.3 in \cite{LOT}, $G_2$ has exactly two vertices $v_1$ and $v_2$.
Assume that $v_1$ and $v_2$ are parallel. Then all the edges of
$G_2$ are positive. Since $M$ is simple, $n_1\geq3$. Thus by Lemma 1.5 in \cite{W3} we can take a label
$x$ which is not a label of a Scharlemann cycle in $G_2$.
Consider the subgraph $\Gamma$ of $G_2$ consisting
of all vertices and all $x$-edges of $G_2$. Let $V,E$ and $F$ be the
numbers of vertices, edges and disk faces of $\Gamma$, respectively.
Since $V<E$, we have $0=\chi(\widehat{F}_2)\le V-E+F < F$, so
$\Gamma$ contains a disk face, which is an $x$-face in $G_2$. This
contradicts \cite[Theorem 4.4]{LOT}.
\end{proof}

We assume without loss of generality that the ordering of the labels
around $v_1$ is anticlockwise while the ordering around $v_2$ is
clockwise. Since the two vertices of $G_2$ are antiparallel, an
Euler characteristic argument on $G_1^+$ shows that there is a disk
face in $G^+_1$, which is a Scharlemann cycle. Then $\widehat{F}_2$
is separating in $M(r_2)$ by Lemma \ref{properties}. We will call
one side of $\widehat{F}_2$ in $M(r_2)$ the $W$(hite) side, denoted
by $\widehat{M}_W$ and the other the $B$(lack) side, denoted by
$\widehat{M}_B$. Then $F_2$ has a $W$ side and a $B$ side in $M$,
which are denoted by $M_W$, $M_B$ respectively. The faces of $G_1$
are called $W$ and $B$ according to which side of $F_2$ a
neighborhood of their boundary lies.

$G_2$ has at most $6$ {\em edge classes}, the isotopy classes in
$\widehat{F}_2$ relative to the fat vertices of $G_2$ by \cite[Lemma 5.2]{G1}. We label
these edge classes $\alpha ,\beta ,\gamma ,\delta ,\epsilon
,\epsilon'$ as in Figure 3. An edge in
$G_1$ or $G_2$ is called an {\em $\eta$-edge} if, being regarded as
an edge in $G_2$, it lies in class $\eta$, $\eta\in\{\alpha ,\beta
,\gamma ,\delta ,\epsilon ,\epsilon'\}$. The $\epsilon$-edges and
$\epsilon'$-edges are positive in $G_2$ while the others are
negative by Lemma \ref{anti}. Since $v_1$ and $v_2$ have the same
valency, $w(\epsilon)=w(\epsilon')$. In Figure 3, an $\epsilon$-edge
has two endpoints on $\partial v_1$. If the endpoint lies between
the edge class $\beta$ and the edge class $\gamma$, we say that the
endpoint lies in class $\epsilon_N$, otherwise $\epsilon_S$. If the
endpoint of an $\epsilon'$-edge lies between the edge class $\alpha$
and the edge class $\gamma$, we say that the endpoint lies in class
$\epsilon'_N$, otherwise $\epsilon'_S$.

\begin{figure}[tbh]
\includegraphics{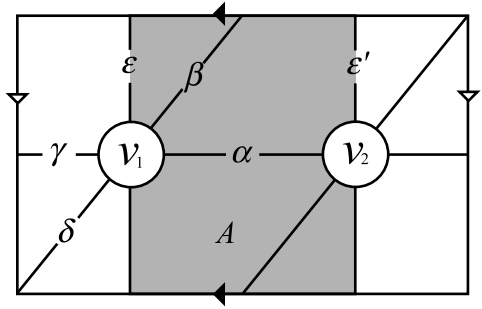}\caption{6 edge classes on $\widehat{F}_2$.}\label{f:jj}
\end{figure}

Two of those (non-loop) edge classes are said to be {\em adjacent}
if ignoring the edge classes $\epsilon, \epsilon'$, the
corresponding edges are adjacent in $\overline{G}_2$. They are {\em
non-adjacent} otherwise. Thus, as in Figure 3, the pair $\{\alpha
,\gamma \}$ and the pair $\{\beta, \delta\}$ are non-adjacent. Let
$\lambda, \mu$ be two elements of $\{\alpha , \beta , \gamma ,
\delta \}$. Then a $(\lambda, \mu )$ {\em face} of $G_1$ is one
whose edges lie in classes $\lambda, \mu$ on $G_2$. Lemma 3.1 of
\cite{GLu2} shows that any such face must have edges in both classes
$\lambda$ and $\mu$. We call any such face a {\em binary face} of
$G_1$.

A $(\lambda, \mu )$ face {\em f} is {\em $\lambda$-good } if no two
consecutive edges in the boundary of {\em f} belong to class
$\lambda$. A $(\lambda, \mu )$ face is {\em good} if it is either
$\lambda$-good or $\mu$-good. A {\em bad} face is a binary face that
is not good. To indicate the coloring of an $(\lambda, \mu )$ face
we will often call it an $X$$(\lambda, \mu )$ face where $X$ is to
take the value $B$ or $W$.

  If $c$ is an $X$ corner of $G_1$ for which the edge incident to
label 1 is in class $\lambda$ and the edge incident to label 2 is in
class $\mu$, we call $c$ an $X$$(\lambda, \mu )$ corner, where $X
\in \{$$B, W$$\}$. Then notice that a bad $X$$(\lambda, \mu )$ face
of $G_1$ must contain $X$$(\lambda, \lambda )$, $X$$(\lambda, \mu
)$, $X$$(\mu, \lambda )$ and $X$$(\mu, \mu )$ corners. For negative
edges, there are four different types of corners i.e. $X(\epsilon_N,
\epsilon'_N)$, $X(\epsilon_N, \epsilon'_S)$, $X(\epsilon_S,
\epsilon'_N)$ and $X(\epsilon_S, \epsilon'_S)$ corners.

Then we have the theorem which plays an important role in this
paper.
\begin{theorem} \label{Se}
Assume that $G_1$ contains a good $X$$(\lambda, \mu )$ face of
length n. Then the $X$ side of $\widehat{F}_2$ in $M(r_2)$ is a
Seifert fibered space over the disk with exactly two exceptional
fibers, one of which has order n. Moreover the Seifert fiber of this
fibration is the curve on $\widehat{F}_2$ formed by $\lambda \cup
\mu$.
\end{theorem}
\begin{proof}
This is exactly Theorem 4.1 in \cite{GL4}.
\end{proof}

The following lemmas are shown in Section 3 of \cite{GL4}.

\begin{lemma}\label{order}
Let $c_1, c_2, \cdot\cdot\cdot, c_n$ be a collection of $X$ corners
of $G_1$, $X\in\{B, W\}$. The anticlockwise ordering of the
endpoints of the edges of these corners on vertex $1$ of $G_2$ is
the same as the clockwise ordering of the edge endpoints on vertex
$2$ of $G_2$.
\end{lemma}

\begin{lemma}\label{same}
$G_1$ cannot contain a vertex $v$ with two edges in the same edge
class incident to $v$ at the same label.
\end{lemma}

\begin{lemma}\label{order1}
Given $X$$\in\{B, W\}$, $G_1$ does not contain $X$$(\lambda, \lambda
)$ corners for three distinct edge classes $\lambda$.
\end{lemma}

\begin{lemma}\label{adjacent}
Suppose that $G_1$ contains a bad $X$$(\alpha, \delta)$ face. Then
there exists $\lambda \in \{\alpha, \delta\}$ such that, for any
edge $e$ of $G_1$ in class $\beta$ or $\gamma$, the edges of $G_1$
adjacent to $e$ on the $X$-side are in class $\lambda$. More
generally, the analogous statement holds whenever
$\{\alpha,\delta\}$ is replaced by any pair of adjacent edge
classes.
\end{lemma}

\begin{lemma}\label{bad}
$G_1$ does not contain bad faces of the same color on distinct edge
class pairs.
\end{lemma}

\begin{lemma}\label{good}
If $G_1$ contains a good $X$$(\lambda, \mu)$ face, then any binary
$X$-face of $G_1$ on an edge class pair distinct but not disjoint
from $\{\lambda,\mu\}$ is bad.
\end{lemma}

\cite[Proposition 4.3]{Lee} says that $G_{1}^{+}$ contains a
connected subgraph $\Lambda$ satisfying the following properties;
\begin{itemize}
\item[(1)] for all vertices $u_x$ of $\Lambda$ but at most one
(an exceptional vertex), there are an edge with edge class $\alpha$
or $\beta$, and an edge with edge class $\gamma$ or $\delta$ in
$\Lambda$ which are incident to $u_x$ with label $i$ for each $i$ =
1,2;
\item[(2)] for the exceptional vertex $u_{x_0}$, if it exists, there are
at least two edges in $\Lambda$ incident to $u_{x_0}$;
\item[(3)] there is a disk $D_{\Lambda}$ $\subset$ $\widehat{F}_1$
such that $IntD_{\Lambda}$ $\cap$ $G_{1}^{+}$ = $\Lambda$; and
\item[(4)] $\Lambda$ has no cut vertex.
\end{itemize}

A vertex of $\Lambda$ is a $boundary$ $vertex$ if there is an arc
connecting it to $\partial D_\Lambda$ whose interior is disjoint
from $\Lambda$, and an $interior$ $vertex$ otherwise. Similarly for
$a$ $boundary$ $edge$ and $an$ $interior$ $edge$ of $\Lambda$.

\begin{lemma} \label{bin}
$\Lambda$ contains a bigon.
\end{lemma}
\begin{proof}
This follows from the proof of \cite[Lemma 4.8]{Lee}.
\end{proof}

\begin{lemma}\label{cy}
$\Lambda$ contains a face bounded only by $\alpha$- or $\beta$-edges
, or only by $\gamma$- or $\delta$-edges. In other words, $\Lambda$
contains a $X$$(\lambda, \mu )$ face where $\{\lambda, \mu \}=
\{\alpha, \beta \}$ or $\{\gamma, \delta \}$ and $X\in \{B, W\}$.
\end{lemma}
\begin{proof}
This follows from \cite[Lemma 4.6]{Lee}.
\end{proof}
By Lemma \ref{cy}, from now on, throughout the paper we may assume
without loss of generality that $\Lambda$ contains a $B(\alpha,
\beta)$ face.

Let $\widehat{H}_B$(resp. $\widehat{H}_W$) be the intersection of
$V_2$ and $\widehat{M}_B$(resp. $\widehat{M}_W$). Let $H_B$(resp.
$H_W$) be the intersection of $V_2$ and $M_B$(resp. $M_W$). Then
$V_2=\widehat{H}_B\cup \widehat{H}_W$ and $\partial V_2=H_B\cup
H_W$. For the purpose of Sections 4 and 5 there are two equivalent
figures of $H_B$ or $H_W$ as in Figure 4 where corners of faces in
$G_1^+$ appear. Figure 4 ($ii$) is obtained by cutting $H_B$ or
$H_W$ along the arc $a$ in Figure 4 ($i$) such that one endpoint of
$a$ lies between the edge class $\alpha$ and the edge class
$\epsilon_S$ at $\partial v_1$ and the other lies between the edge
class $\alpha$ and the edge class $\beta$ at $\partial v_2$.

\begin{figure}[tbh]
\includegraphics{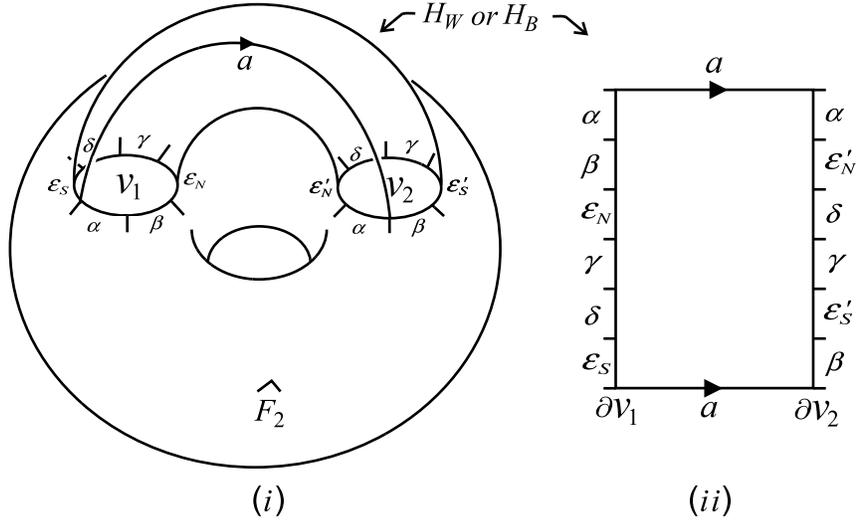}\caption{Two equivalent figures of $H_W$ or $H_B$.}\label{f:jj}
\end{figure}

\begin{lemma} \label{abface}
A $B(\alpha, \beta)$ face is good.
\end{lemma}
\begin{proof}
Suppose a $B(\alpha, \beta)$ face is bad. Then it contains $B(\alpha, \alpha)$,
$B(\alpha, \beta)$, $B(\beta, \alpha)$, $B(\beta, \beta)$ corners.
These corners lie on $H_B$ as in Figure 5. There are two cases.
First, we assume Figure 5 ($i$). Figure 5 ($ii$) is similar. It follows from Figure 5 ($i$) that
$\omega(\epsilon')+ \omega(\delta)+ \omega(\gamma)+ \omega(\epsilon') \leq \omega(\beta)$.
On the other hand, $\omega(\alpha)+ \omega(\epsilon')+ \omega(\delta)+ \omega(\gamma)+ \omega(\epsilon')+
\omega(\beta)= 3n_1$ since the number of all the labels around a vertex in $G_2$ is 3$n_1$.
Therefore, combining this equation with the above inequality, we get the inequality $\omega(\alpha)+ 2\omega(\beta) \geq 3n_1$.
This implies that either $\omega(\alpha) \geq n_1$ or $\omega(\beta) \geq n_1$, which contradicts \cite[Lemma 1.5]{W3}.
\end{proof}

\begin{figure}[tbh]
\includegraphics{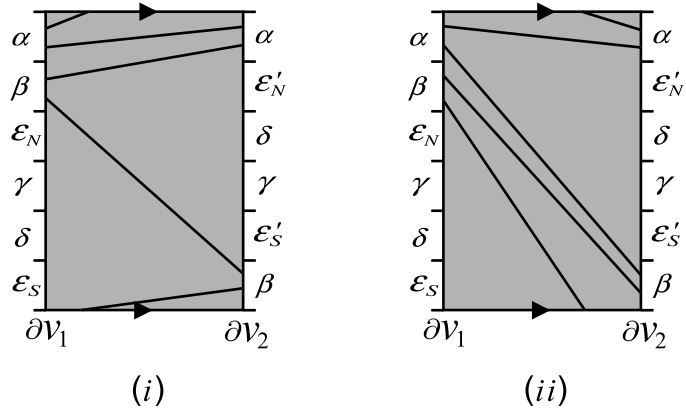}\caption{The intersection of $H_B$ and the corners of a $B(\alpha, \beta)$ face.}\label{f:jj}
\end{figure}
\begin{lemma} \label{Se1}
The black side $\widehat{M}_B$ of $M(r_2)$ is a Seifert fibered
space over the disk with exactly two exceptional fibers. Moreover
the Seifert fiber of this fibration is the curve on $\widehat{F}_2$
formed by $\alpha \cup \beta$.
\end{lemma}
\begin{proof}
$\Lambda$ contains a $B(\alpha, \beta)$ face, which we denote by
$f$. Let $\lambda_i$, $i=$ 1,$\cdot\cdot\cdot$, $n$,  be the corner
of $f$ and $\partial\lambda_{i}^{j}$ the endpoint of $\lambda_i$ at
the label $j$, $j=$ 1, 2. Let $I(v_j)$ be the shortest interval on
$\partial v_j$ containing all $\partial \lambda_{i}^{j}$ in its
interior. Relabel $i$ so that $\partial I(v_1)=$ $\{\partial
\lambda_{1}^{1},\partial \lambda_{n}^{1}\}$. It follows from
\cite[Lemma 2.4]{G1} that $\partial I(v_2)=$ $\{\partial
\lambda_{1}^{2},\partial \lambda_{n}^{2}\}$. Since $f$ is good by Lemma \ref{abface}, there exists a
disk $D$ in $H_B$ bounded by four arcs $a_1$, $a_2$, $b_1$, $b_2$
where $a_i$ is the arc of $\partial v_i$ containing $I(v_i)$, $i= 1,
2$, and $b_j$ is the parallel arc to $\lambda_j$ in $H_B$, $j= 1,2$,
so that $D$ contains all the corners of $f$.

We can take an annulus $A$ on $F_2$ containing all the edges of $f$
so that $A\cup D$ is an once punctured torus on $M_B$ and the two
boundary components of $A$ are homotopic to the curve on
$\widehat{F}_2$ formed by $\alpha \cup \beta$. See Figure 3. Since
$\partial f$ is non-separating on $A\cup D$, we can get a disk $D'$
by surgering $A\cup D$ along $f$. Push $D'$ slightly into the
interior of $M_B$ so that $D'$ is a properly embedded disk in $M_B$.
Notice that $\partial D'$ contains the two arcs $b_1$, $b_2$. An
orientation of $\partial D'$ induces orientations of $b_1$, $b_2$
that are opposite on $H_B$. Hence shrinking $\widehat {H}_B$ to its
core yields a properly embedded annulus $A'$ in the black side
$\widehat{M}_B$ of $\widehat{F}_2$. In addition, $A'$ separates
$\widehat{M}_B$ into two manifolds $M_1$ and $M_2$, both of which
have a torus boundary i.e. $\widehat{M}_B$ = $M_1 \cup_{A'} M_2$. We
can observe that the core of $\widehat{H}_B$ lies on $A'$. Therefore
by isotoping the core of $V_2$ slightly off $\partial M_1$ and
$\partial M_2$, we can show that $M_1$ and $M_2$ are solid tori by
the choice of $\widehat{F}_2$. Since $\widehat{F}_2$ is essential in
$M(r_2)$, the black side $\widehat{M}_B$ is a Seifert fibered space
over the disk with two exceptional fibers and the fiber is the core
of $A'$ and hence the core of $A$, as desired.
\end{proof}

$Remark$: In the proof of Lemma \ref{Se1}, the core of
$\widehat{H}_B$, that is, the black part of $V_2$, lies on the
annulus which separates the Seifert fibered space $\widehat{M}_B$
into two solid tori. We will use this fact to prove Lemma
\ref{whit}.

Before proving the next lemma, we observe that the edge classes of
the three edges of a trigon in $\Lambda$ cannot be all distinct by
Lemma \ref{order} and cannot be all same by Lemma 3.1 of \cite{GLu2}.
\begin{lemma} \label{whit}
If $\Lambda$ contains a white bigon or a white trigon, the two edge
classes of the edges of a white bigon or a white trigon are
$\lambda$ and $\mu$ where $\{\lambda, \mu \}$ $=$ $\{\alpha, \delta \}$, $\{\alpha, \gamma \}$,
$\{\beta, \delta \}$ or $\{\beta, \gamma \}$.
\end{lemma}
\begin{proof}
Let $f$ be a white bigon in $\Lambda$. $f$ is necessarily a good
$W$$(\lambda, \mu)$ face. Suppose $\{\lambda, \mu \}$ = $\{\alpha,
\beta \}$ or $\{\gamma, \delta \}$. Then applying the same argument
as in the proof of Lemma \ref{Se1} using $f$, we see that the
white side $\widehat{M}_W$ is a Seifert fibered space over the disk
with exactly two exceptional fibers and the Seifert fiber of this
fibration is the curve on $\widehat{F}_2$ formed by $\lambda \cup
\mu$. In other words, $\widehat{M}_W=M_3 \cup_{A} M_4$ where $M_3$
and $M_4$ are solid tori and $A$ is a properly embedded annulus in
$\widehat{M}_W$, and the core of $A$ is the fiber. On the
other hand, the black side $\widehat{M}_B=M_1 \cup_{A'} M_2$ where
$M_1$ and $M_2$ are solid tori and the core of $A'$ is the fiber. By
the $Remark$ below Lemma \ref{Se1}, the core of $V_2$ lies on $A\cup
A'$. Furthermore since two fibers of $\widehat{M}_B$ and
$\widehat{M}_W$ coincide on $\widehat{F}_2$, we can isotope $A$ and
$A'$ on $\widehat{F}_2$ fixing the core of $V_2$ so that $A\cup A'$
is a torus $\widehat{T}$, say. Note that either side of
$\widehat{T}$ is a Seifert fibered space over the disk with two
exceptional fibers, which implies that $\widehat{T}$ is essential in
$M(r_2)$. Pushing the core of $V_2$ off $\widehat{T}$, we get an
essential torus in $M(r_2)$ which misses the core of $V_2$, which is
a contradiction to the choice of $\widehat{F}_2$. Therefore
$\{\lambda, \mu \}$ = $\{\alpha, \delta \}$, $\{\alpha, \gamma \}$,
$\{\beta, \delta \}$ or $\{\beta, \gamma \}$.

Let $g$ be a white trigon in $\Lambda$. Hence $g$ is also
necessarily a good $W$$(\lambda, \mu)$ face. Applying the same
argument as in the case of a bigon above, we complete the proof.
\end{proof}

\begin{lemma} \label{twobi}
Any two black bigons in $G_{1}^{+}$ have the same pair of edge
classes. Furthermore the pair is either $\{\alpha, \beta \}$ or
$\{\gamma, \delta \}$.
\end{lemma}
\begin{proof}
Lemma \ref{order} implies that no two black bigons of $G_{1}^{+}$
have distinct and non-disjoint pairs of edge classes, nor have the
pair $\{\{\alpha, \gamma\}, \{\beta, \delta\}\}$. Therefore the only
possible different pairs of edge classes of two bigons are
$\{\{\alpha, \beta\}, \{\gamma, \delta\}\}$, $\{\{\alpha, \delta\},
\{\beta, \gamma\}\}$. In either case, we construct two M$\ddot{\rm{o}}$bius bands
from the two bigons by shrinking $\widehat{H}_B$ to its core.
Furthermore, the two boundaries of the two M$\ddot{\rm{o}}$bius bands cobound an
annulus in $\widehat{F}_2$, which gives a Klein bottle in $M(r_2)$,
a contradiction to the assumption that $M(r_2)$ does not contain a Klein
bottle. Hence any two black bigons in $G_{1}^{+}$ have the same
pair of edge classes.

Now we show that the pair is either $\{\alpha, \beta \}$ or
$\{\gamma, \delta \}$. If there is a black bigon in $G_{1}^{+}$ with
the edge class pair $\{\lambda, \mu \}$, then the black side
$\widehat{M}_B$ of $\widehat{F}_2$ in $M(r_2)$ is a Seifert fibered
space over the disk with exactly two exceptional fibers and the
Seifert fiber of this fibration is the curve on $\widehat{F}_2$
formed by $\lambda \cup \mu$ by Theorem \ref{Se}. However the fiber
of $\widehat{M}_B$ is the curve on $\widehat{F}_2$ formed by $\alpha
\cup \beta$ by Lemma \ref{Se1}. Since $M(r_2)$ does not contain a Klein
bottle, the Seifert fibration of $\widehat{M}_B$ is unique by \cite
[Theorem VI.18]{J}. Therefore the edge class pair $\{\lambda, \mu
\}$ of the bigon must be either $\{\alpha, \beta \}$ or $\{\gamma,
\delta \}$.
\end{proof}

The next theorem is the main result of the combinatorial arguments
in this paper. To prove this, we use Lemma \ref{trigon1} and the main theorems in Sections 4
and 5.
\begin{theorem} \label{whitebigon}
$\Lambda$ contains a $W(\lambda, \mu)$ bigon, where $\{\lambda, \mu
\}=\{\alpha, \delta \}$, $\{\alpha, \gamma \}$, $\{\beta, \delta \}$
or $\{\beta, \gamma \}$.
\end{theorem}
\begin{proof}
Lemma \ref{bin} says that $\Lambda$ contains a bigon $f$. If $f$ is
white, then we are done by Lemma \ref{whit}.

Suppose $\Lambda$ does not contain a white bigon. Then $f$ is black. By
Lemma \ref{twobi}, we may assume that $f$ is a $B(\alpha, \beta)$
bigon. Lemma \ref{trigon1} says that $\Lambda$ contains a trigon. On
the other hand, Theorem \ref{noblacktrigon} and Theorem
\ref{nowhitetrigon} imply that $\Lambda$ does not contain a trigon. This
contradiction shows that $\Lambda$ contains a white bigon as
desired.
\end{proof}

To prove Lemma \ref{trigon1}, Theorem \ref{noblacktrigon} and
Theorem \ref{nowhitetrigon}, we assume that $\Lambda$ contains a
$B(\alpha, \beta)$ bigon and does not contain a white bigon.

Lemma \ref{twobi} enables us to apply Lemmas 4.10 and 4.11 of
\cite{Lee} without change.
\begin{lemma}\label{valence}
All interior vertices of $\Lambda$ have valency at least $4$ in
$\overline{\Lambda}$ and all boundary vertices of $\Lambda$ but the
exceptional vertex have valency at least $3$ in
$\overline{\Lambda}$.
\end{lemma}
\begin{proof}
This is Lemma 4.10 of \cite{Lee}.
\end{proof}

\begin{lemma}\label{trigon1}
$\Lambda$ contains a trigon.
\end{lemma}
\begin{proof}
This is Lemma 4.11 of \cite{Lee}.
\end{proof}

\begin{lemma} \label{interior}
$\Lambda$ contains an interior vertex.
\end{lemma}
\begin{proof}
Suppose that $\Lambda$ has no interior vertices. Lemma \ref{valence}
guarantees that the number of vertices of $\Lambda$ is greater than
$3$. Put $\Lambda$ on an abstract disk $D$ so that $\Lambda$ lies in
the interior of $D$ and insert an edge connecting each vertex of
$\Lambda$ to $\partial D$. By taking the reduced graph of the
resulting graph, we obtain a graph $\Gamma$ satisfying the
supposition of \cite[Lemma 2.6.5]{CGLS}. In addition, since the
number of vertices of $\Gamma$ is greater than $3$, $\Gamma$ must
satisfy (*) in the proof of \cite[Lemma 2.6.5]{CGLS}. In other
words, $\Gamma$ contains two vertices of valency at most 3. However
this is impossible since all vertices of $\Gamma$, possibly except
one, have valency at least $4$ by Lemma \ref{valence}.
\end{proof}

In Section 4, we will show that $\Lambda$ does not contain a black trigon
and in Section 5, we will show that $\Lambda$ does not contain a white
trigon.

\section{a $B(\alpha, \beta)$ bigon and a black trigon}

Throughout this section, we assume that $\Lambda$ contains a
$B(\alpha, \beta)$ bigon and does not contain a white bigon. Our goal in
this section is to show the following.

\begin{theorem} \label{noblacktrigon}
If $\Lambda$ contains a $B(\alpha, \beta)$ bigon and does not contain a
white bigon, then $\Lambda$ does not contain a black trigon.
\end{theorem}

The main idea is to define two dual graphs of $\Lambda$ and use the
fact that there is no triple of non-isomorphic disk faces in
$\Lambda$ and apply the index equation $\sum _{vertices} I(v)$ +
$\sum_{faces} I(f)=2$.

To prove Theorem \ref{noblacktrigon}, we suppose throughout this section that $\Lambda$ contains a black trigon.

\begin{lemma} \label{trigons}
All black trigons of $G_{1}^{+}$ have the same edge classes, i.e.
two $\gamma$-edges and one $\delta$-edge, say.
\end{lemma}
\begin{proof}
Let $g$ be a black trigon of $G_{1}^{+}$. Lemma \ref{order} excludes
the case that the three edges of $g$ have distinct edge classes.
Thus $g$ has two $\lambda$-edges and one $\mu$-edge. Since $g$ is a
good $B(\lambda, \mu)$ face, we have the Seifert fibration of the
black side $\widehat{M}_B$ resulting from $g$ whose fiber is
represented by $\{\lambda, \mu \}$. Since $M(r_2)$ does not contain a Klein
bottle, by the uniqueness of the Seifert fibration, $\{\lambda,
\mu \}=\{\alpha, \beta \}$ or $\{\gamma, \delta \}$ by Lemma
\ref{Se1}.

Since $\Lambda$ contains a $B(\alpha, \beta)$ bigon,
there are two corners, a $B(\alpha,\beta)$ corner and
a $B(\beta, \alpha)$ corner. Then Lemma \ref{order} eliminates the case that
$\{\lambda, \mu \}$ = $\{\alpha, \beta \}$. However Lemma
\ref{order} also excludes the case that $\Lambda$ contains
simultaneously two trigons with two $\gamma$-edges and one
$\delta$-edge, and two $\delta$-edges and one $\gamma$-edge
respectively. Hence we may assume that all trigons of $G_{1}^{+}$
have the same edge classes, i.e. without loss of generality two
$\gamma$-edges and one $\delta$-edge.
\end{proof}

\begin{lemma} \label{elimination} $\Lambda$ cannot contain a
$B(\lambda, \mu)$ corner, where $(\lambda, \mu)= (\alpha, \delta)$,
$(\beta, \gamma)$, $(\gamma, \beta)$, $(\delta, \alpha)$, or
$(\delta, \delta)$.
\end{lemma}
\begin{proof}
$\Lambda$ contains a black bigon $f_1$ with one $\alpha$-edge and
one $\beta$-edge, and contains a black trigon $f_2$ with two
$\gamma$-edges and one $\delta$-edge. Consider the intersection of
all the corners of $f_1$ and $f_2$, and $H_B$. See Figure 6. Then
the above black corners are impossible.
\end{proof}

\begin{figure}[tbh]
\includegraphics{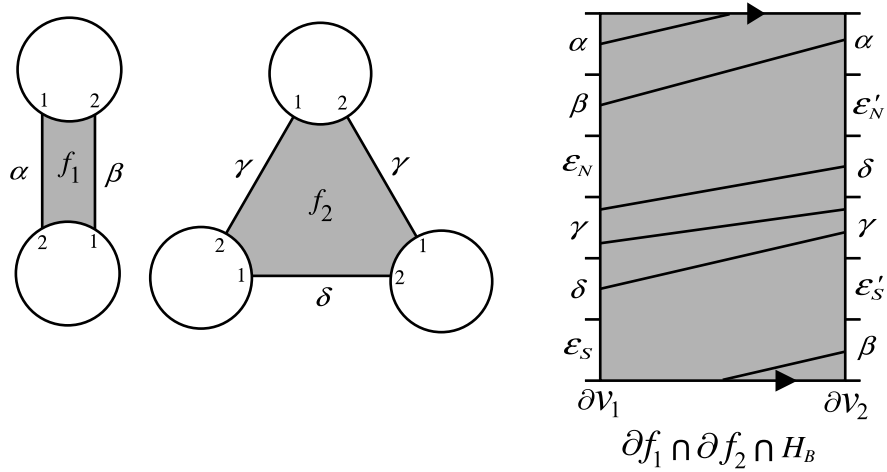}\caption{A $B(\alpha, \beta)$ bigon, a $B(\gamma,
\delta)$ trigon and the intersection on $H_B$.}\label{f:jj}
\end{figure}

Since $\Lambda$ contains a black face and a white face
simultaneously by Lemma \ref{valence}, $M_B$ and $M_W$ are
handlebodies of genus 2 from \cite [Lemma 8.3] {GLu2}.

Traveling around the boundary of a disk face of $G_1^+$ gives rise
to a cyclic sequence of edge classes. We shall say that two disk
faces of $G_1^+$ of the same color are {\em isomorphic} if the
cyclic sequences obtained by traveling in some directions are equal.
See \cite[Section 5]{GLu4}.

Recall that $H_B=V_2$ $\cap$ $M_B$ and $H_W=V_2$ $\cap$ $M_W$. Then
$\partial M_B=H_B$ $\cup$ $F_2$ and $\partial M_W=H_W$ $\cup$ $F_2$.
Since all the vertices of $\Lambda$ are parallel, each face of
$\Lambda$ is a non-separating disk in $M_B$ or $M_W$. Note that any
two faces of the same color in $\Lambda$ are disjoint.

\begin{lemma}\label{word}
If two disk faces of $G_1^+$ are parallel in $M_B$ or $M_W$, then
they are isomorphic.
\end{lemma}
\begin{proof}
This is Lemma 4.13 of \cite{Lee}.
\end{proof}

\begin{theorem} \label{notriple}
There is no triple of mutually non-isomorphic black disk faces in
$G_{1}^+$, and there is no triple of mutually non-isomorphic white
disk faces in $G_{1}^+$.
\end{theorem}

To prove this theorem, we will show that there is a white binary
disk face by using dual graphs of $\Lambda$.

For the purpose of the proof of Theorem \ref{notriple}, we assume
for the time being that $G_1= G_1^+$. Then every vertex of $\Lambda$
except an exceptional vertex has valency 6 since all edges of $G_1$
are positive.

Now we define a dual graph $\Lambda^*$ to $\Lambda$ as in
\cite[Section 3]{GL4}. To do this, it is convenient to regard
$\Lambda$ as a graph in $S^2$, rather than a disk. In other words,
we add an additional $outside$ $face$ to the set of faces of
$\Lambda$. Then a $vertex$ of $\Lambda^*$ is a point in the interior
of each face of $\Lambda$ and a point of the outside face. An $edge$
of $\Lambda^*$ is the edge connecting the vertices of $\Lambda^*$
corresponding to the faces of $\Lambda$ on either side of each edge
of $\Lambda$, and meeting each edge transversely in a single point.

We give two orientations $\omega$ and $\omega'$ to the edges of
$\Lambda^*$ as follows;
\begin{itemize}
\item[$\omega:$] \em{an edge of $\Lambda^*$ dual to an edge of $\Lambda$
in an edge class $\in\{\alpha, \delta \}$ $($resp. $\{\beta, \gamma
\}$$)$ is oriented from the $W$-side to the $B$-side $($resp. from
the $B$-side to the $W$-side$)$}.
\end{itemize}
\begin{itemize}
\item[$\omega':$] \em{an edge of $\Lambda^*$ dual to an edge of $\Lambda$
in an edge class $\in\{\beta, \delta \}$ $($resp. $\{\alpha, \gamma
\}$$)$ is oriented from the $W$-side to the $B$-side $($resp. from
the $B$-side to the $W$-side$)$}.
\end{itemize}
See Figure 7.

 $Note$. For an edge $e$ of $\Lambda$ in the boundary
of the outside face, we regard the side of $e$ contained in the
outside face as being locally colored with the color opposite to
that of the face of $\Lambda$ that contains $e$.

\begin{figure}[tbh]
\includegraphics{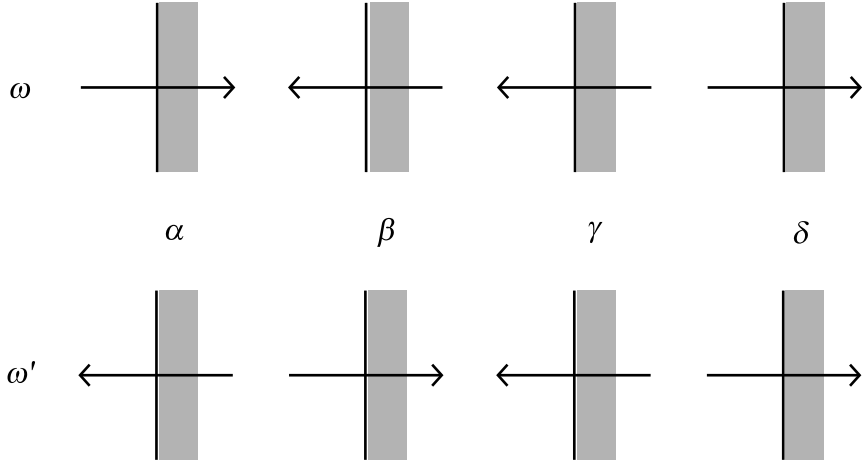}\caption{Two orientations $\omega$ and $\omega'$ to the edges of $\Lambda^*$.}\label{f:jj}
\end{figure}

As defined in \cite{GLu1}, the $index$ of a vertex $v$ and a face
$f$ of a graph in $S^2$ with oriented edges is defined by
$I(v)=1-s(v)/2$ and $I(f)=\chi(f)-s(f)/2$ where $s(v)$ and $s(f)$
are switches at $v$ and of $f$ respectively. A $sink$ or a $source$
is a vertex of index 1, and a $cycle$ is a face of index 1. We say
an $exceptional$ $cycle$ if it is dual to an exceptional vertex.

\begin{lemma} \label{nocycle}
$\Lambda^*$ does not contain a cycle except an exceptional cycle with
respect to $\omega, \omega'$.
\end{lemma}
\begin{proof}
Suppose $v$ is a non-exceptional vertex dual to a cycle with respect
to $\omega$. Since $v$ has valency 6, the edge classes of three
edges which are incident to $v$ at the label 1 belong to either
$\{\alpha, \delta \}$ or $\{\beta, \gamma \}$. Therefore two of them
have the same edge class at the same label. This contradicts Lemma
\ref{same}. The same argument applies to $\omega'$.
\end{proof}

Let $s$ and $s'$ be the number of sinks and sources of $\Lambda^*$
at vertices dual to the faces of $\Lambda$ (i.e. we exclude the
outside face) with respect to $\omega, \omega'$ respectively.

\begin{lemma} \label{ss}
$s, s'$ $>$ $0$.
\end{lemma}
\begin{proof}
By Lemma \ref{nocycle}, the only possible cycle is an exceptional
cycle. However if $\Lambda^*$ has an exceptional cycle, then it cannot
contain the sink or source dual to the outside face at the same
time. Hence it follows from the equation $\sum _{vertices} I(v)$ +
$\sum_{faces} I(f)=2$ that $s, s'$ $>$ $0$.
\end{proof}

A sink or source of $\omega, \omega'$ in $\Lambda^*$ that does not correspond to
the outside face of $\Lambda$ is a binary face of $\Lambda$.

\begin{lemma} \label{whitebinary}
There is a white binary face in $\Lambda$ dual to a sink or source with respect to $\omega$.
\end{lemma}
\begin{proof}
By Lemma \ref{ss}, there is a binary face in $\Lambda$ dual to a sink or source
with respect to $\omega$. We show that this binary face is white.

Suppose that a binary face in $\Lambda$ dual to a sink or source with
respect to $\omega$ is black. Then such a black binary face is either a $B(\alpha, \delta)$ face
or a $B(\beta, \gamma)$ face. The former contains a $B(\alpha, \delta)$ corner and
the latter contains a $B(\beta, \gamma)$ corner. This contradicts Lemma \ref{elimination}.
This completes the proof.
\end{proof}

So far, we have shown that under the assumption that $G_1=G^+_1$ Lemma \ref{whitebinary}
holds. Now we are ready to prove Theorem \ref{notriple}.
\begin{proof}[Proof of Theorem \ref{notriple}]
We follow the proof of Lemma 4.16 of \cite{Lee}. Let $f_1, f_2, f_3$
be three mutually non-isomorphic black faces of $G_{1}^{+}$. Then
these faces cut $M_B$ into two $3$-balls by Lemma \ref{word}.

 Let $f$ be a black face of $G_1$ other than $f_1,f_2,f_3$. Since $f$
lies in the complement of $f_1$ $\cup$ $f_2$ $\cup$ $f_3$ in $M_B$,
$f$ must be a disk, otherwise it would be compressible in $M_B$, so
$F_1$ would be compressible in $M$, a contradiction to the fact that
$F_1$ is incompressible. Each component of $\partial f \cap F_2$ is
an edge of $G_2$. $\partial f$ must be an essential curve in
$\partial M_B=F_2\cup H_B$, otherwise some component of $\partial f
\cap F_2$ would be a trivial loop in $G_2$. Hence $f$ is an
essential disk in $M_B$, so it must be parallel to one of the faces
$f_1,f_2$ and $f_3$. Therefore any black face of $G_1$ is a disk
face isomorphic to one of $f_1,f_2,f_3$ by Lemma \ref{word}. It
follows that all the edges of $G_1$ are positive, i.e. $G_1=G_1^+$.

Next, we will show that any white face of $G_1$ is a disk face. Note
that if every white face of $G_1$ is a disk face, then $G_1$ is
connected.  Suppose a white face $g$ of $G_1$ is not a disk face.
Lemma \ref{whitebinary} guarantees the existence of a white binary
disk face $g_1$. Furthermore, there must be a white disk face $g_2$
other than $g_1$ in $\Lambda$ which is non-isomorphic to $g_1$,
otherwise there are three edges in the same edge class incident to
an interior vertex, contradicting Lemma \ref{same}. Thus the complement of $g_1\cup g_2$ in $M_W$ is
a 3-ball. So, $g$ is compressible in $M_W$, which implies that $F_1$
is compressible, a contradiction.

We have shown that $G_1=G_1^+$ and $G_1$ is connected. Therefore,
$F_1$ is nonseparating, which contradicts \cite[Theorem 1.1]{Lee2}.
This completes the proof that there is no triple of mutually
non-isomorphic black disk faces in $G_{1}^+$. Similarly, we can
apply the same argument to show that there is no triple of
non-isomorphic white disk faces, by using a $B(\alpha, \beta)$
bigon.
\end{proof}
Now we go back to the subgraph $\Lambda$ as in Section 3 and
continue with the proof of Theorem \ref{noblacktrigon}.
Since $\Lambda$ contains a $B(\alpha, \beta)$ bigon and a $B(\gamma,
\delta)$ trigon with two $\gamma$-edges and one $\delta$-edge by
Lemma \ref{trigons}, Theorem \ref{notriple} implies the following.

\begin{lemma} \label{ast}
Any black disk face of $G_{1}^{+}$ is isomorphic to either a
$B(\alpha, \beta)$ bigon or a $B(\gamma, \delta)$ trigon with two
$\gamma$-edges and one $\delta$-edge.
\end{lemma}

We define
the dual graph $\Lambda^*$ with two orientations $\omega, \omega'$
as in the case $G_1 = G_1^+$. The cycle dual to an interior
vertex in $\Lambda$ is said to be an $ordinary$ $cycle$. The cycle
dual to a boundary vertex (except an exceptional vertex) is called a
$boundary$ $cycle$. The cycle dual to an exceptional vertex is said
to be an $exceptional$ $cycle$. Thus when we say a boundary cycle
dual to a boundary vertex, we do not allow it to be the cycle dual
to an exceptional vertex.

\begin{lemma} \label{noblackss}
$\Lambda^*$ does not contain a sink or source dual to a black face of
$\Lambda$ with respect to $\omega, \omega'$.
\end{lemma}
\begin{proof}
Lemma \ref{ast} says that any black disk face is either a $B(\alpha,
\beta)$ bigon or a $B(\gamma, \delta)$ trigon with two
$\gamma$-edges and one $\delta$-edge. However, either black disk
face is not dual to a sink or source with respect to $\omega,
\omega'$.
\end{proof}

It is guaranteed by Lemma \ref{interior} that $\Lambda$ contains an
interior vertex.

Note that there are only four distinct types of interior vertices.
See Figure 8.

\begin{figure}[tbh]
\includegraphics{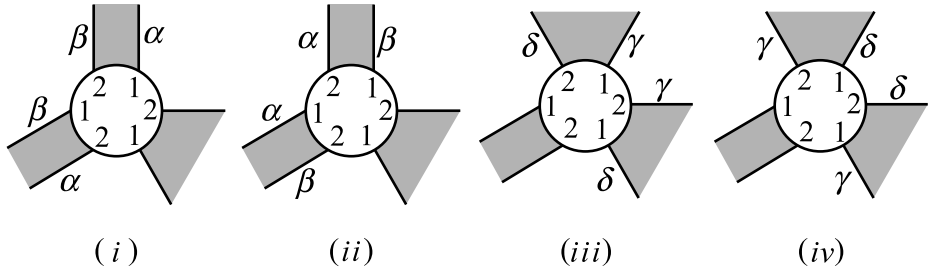}\caption{Four distinct types of interior vertices.}\label{f:jj}
\end{figure}

\begin{lemma} \label{noordi}
$\Lambda^*$ does not contain an ordinary cycle with respect to
$\omega,\omega'$.
\end{lemma}
\begin{proof}
Since any interior vertex has valency 6, this follows from the proof
of Lemma \ref{nocycle}.
\end{proof}

\begin{lemma} \label{bdrycycle}
A boundary cycle dual to a boundary vertex with respect to $\omega,
\omega'$ is one of those illustrated in Figure $\mathrm{9}$ and
Figure $10$ respectively.
\end{lemma}

\begin{proof}
If a boundary vertex $v$ dual to a boundary cycle has valency 5 or
6, then at least two edges incident to $v$ have the same edge class
at the same label, contradicting Lemma \ref{same}. Hence any
boundary vertex dual to a boundary cycle has valency 4. Furthermore,
any boundary vertex has two $\alpha$- or $\beta$-edges and also two
$\gamma$- or $\delta$-edges by the property (1) of $\Lambda$. Then
the proof follows from Lemma \ref{ast}.
\end{proof}
\begin{figure}[tbh]
\includegraphics{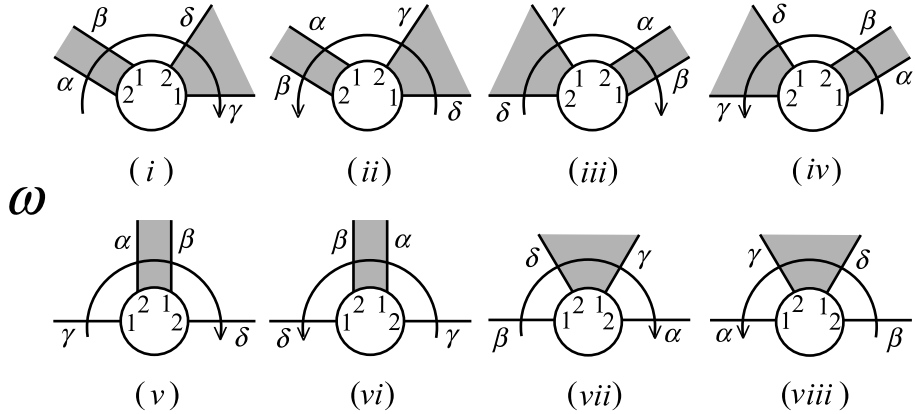}\caption{Boundary cycles with respect to $\omega$.}\label{f:jj}
\end{figure}
\begin{figure}[tbh]
\includegraphics{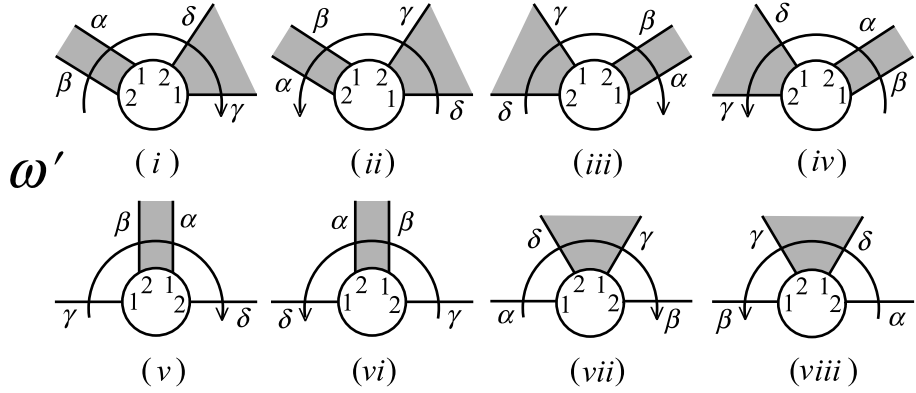}\caption{Boundary cycles with respect to $\omega'$.}\label{f:jj}
\end{figure}

\begin{lemma} \label{omega2}
$\Lambda^*$ does not contain a boundary cycle with respect to $\omega'$.
\end{lemma}
\begin{proof}
Assume for contradiction that $\Lambda$ contains a boundary cycle
dual to a boundary vertex $u_x$. In $G_2$, the label $x$ appears
either three times around $v_1$ at the ends of $\alpha$-, $\gamma$-
and $\epsilon$-edges and three times around $v_2$ at the ends of
$\beta$-, $\delta$- and $\epsilon'$-edges, or three times around
 $v_1$ at the ends of $\beta$-, $\delta$- and $\epsilon$-edges and
three times around $v_2$ at the ends of $\alpha$-, $\gamma$- and
$\epsilon'$-edges. The odd-numbered types in Figure 10 have the first
case and the even-numbered types have the second case.

\begin{claim}
In $\mathrm{G}_2$, the label $x$ appears at $\epsilon_N$ (resp. $\epsilon_S$) of the
$\epsilon$-edge class at $v_1$ if and only if it appears at
$\epsilon'_S$ (resp. $\epsilon'_N$) of the $\epsilon'$-edge class at $v_2$.
\end{claim}
\begin{proof}
Assume the first case i.e. the label $x$ appears three times around
the vertex $v_1$ at the ends of $\alpha$-, $\gamma$- and
$\epsilon$-edges and three times around the vertex $v_2$ at the ends
of $\beta$-, $\delta$- and $\epsilon'$-edges. The second case is
similar.

\begin{figure}[tbh]
\includegraphics{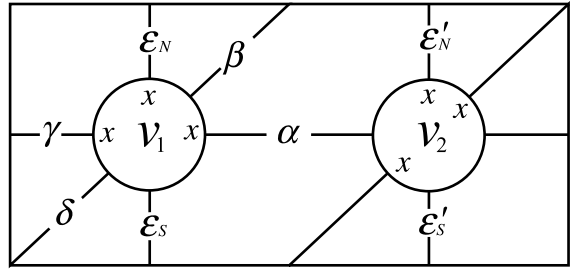}\caption{Three occurrences of label $x$ on $v_1$ and $v_2$.}\label{f:jj}
\end{figure}

Suppose, for example, that $x$ appears at $\epsilon_N$ and at
$\epsilon'_N$ of $\epsilon$- and $\epsilon'$-edges respectively at
the same time. See Figure 11. At the vertex $v_1$ the edge classes
$\delta$ and $\epsilon$ don't have the label $x$. Then the sum of
the weights of $\delta$ and $\epsilon$ is less than $n_1$. On the
other hand, at the vertex $v_2$ the edge classes $\delta$ and
$\epsilon'$ have the label $x$ twice, which implies that the sum of
the weights of $\delta$ and $\epsilon'$ is greater than $n_1$.
However, these two inequalities conflict since the weights of
$\epsilon$ and $\epsilon'$ are same.
\end{proof}

Assume that the boundary vertex $u_x$ is dual to one of the types
($i$)$\sim$($iv$). Then we have a $B(\epsilon,
\epsilon')$ corner. In Figure 6 of Lemma \ref{elimination},
($\epsilon,\epsilon'$) should be ($\epsilon_N, \epsilon'_N$) or
($\epsilon_S, \epsilon'_S$). This contradicts the Claim.

For the type ($v$), we have two black corners, a $B(\gamma,
\epsilon')$ corner, a $B(\epsilon, \delta)$ corner. However, from
Figure 6 of Lemma \ref{elimination}, ($\epsilon, \epsilon'$) should
be ($\epsilon_N, \epsilon'_N$), which is a contradiction to the
Claim. The other types are also impossible by the same argument.
This completes the proof.
\end{proof}

$Note$. The Claim in the proof of Lemma \ref{omega2} doesn't work
with the orientation $\omega$.

\begin{lemma} \label{omega1}
$\Lambda^*$ does not contain a boundary cycle with respect to $\omega$.
\end{lemma}

\begin{proof}
Lemmas \ref{noblackss}, \ref{noordi}, \ref{omega2} and the equation
$\sum _{vertices} I(v)$ + $\sum_{faces}$ $I(f)=2$ imply that there
is a sink or source dual to a white disk face with respect to
$\omega'$. In other words, there is a $W(\beta, \delta)$ or
$W(\alpha, \gamma)$ face. There are eight types of boundary cycles
with respect to $\omega$ as in Figure 9. The first four types have a
$B(\epsilon, \epsilon')$ corner. However, by Figure 6,
($\epsilon, \epsilon'$) must be either ($\epsilon_N, \epsilon'_N$)
or ($\epsilon_S, \epsilon'_S$). Assume we have the type ($i$).
If we let $u_x$ be a boundary vertex of the type $(i)$, then the
label $x$ appears three times around the vertex $v_1$ at the ends of
$\beta$-, $\gamma$- and $\epsilon$-edges and three times around the
vertex $v_2$ at the ends of $\alpha$-, $\delta$- and
$\epsilon'$-edge. Applying a similar argument as in the Claim in
Lemma \ref{omega2}, ($\epsilon, \epsilon'$) must be ($\epsilon_S,
\epsilon'_S$). Thus, we have a $W(\beta, \delta)$ corner, a
$W(\epsilon_S,\alpha)$ corner, and a $W(\gamma, \epsilon'_S)$
corner. Since we have a $W(\beta, \delta)$ or $W(\alpha, \gamma)$
face, there are $W(\delta, \delta)$ or $W(\beta, \beta)$
corners in a $W(\beta, \delta)$ face by the Remark below the proof of
Lemma 4.4 in \cite{GL4}, and there are $W(\alpha, \gamma)$
corners in a $W(\alpha, \gamma)$ face. Then neither a $W(\beta,
\beta)$ corner, a $W(\delta, \delta)$ corner nor a $W(\alpha,
\gamma)$ corner satisfies Lemma \ref{order} with $W(\beta,
\delta)$, $W(\epsilon_S,\alpha)$, and $W(\gamma, \epsilon'_S)$
corners. Therefore type ($i$) can't happen. A similar argument
applies to the other types ($ii$), ($iii$) and ($iv$).

For the type ($v$), we have a $W(\beta, \delta)$ corner, a
$W(\gamma, \alpha)$ corner and a $W(\epsilon, \epsilon')$ corner.
Then Lemma \ref{order} forces ($\epsilon, \epsilon'$) to be either
($\epsilon_N, \epsilon'_S$) or ($\epsilon_N, \epsilon'_S$). If
($\epsilon, \epsilon'$) = ($\epsilon_N, \epsilon'_S$), the boundary
vertex has a $B(\gamma, \epsilon'_S)$ corner, which is impossible by
Figure 6. If ($\epsilon, \epsilon'$) = ($\epsilon_S, \epsilon'_N$),
then the boundary vertex has a $B(\epsilon_S, \delta)$ corner,
which is also impossible by Figure 6. For the other types, we can
apply a similar argument.
\end{proof}

\begin{proof}[Proof of Theorem $\ref{noblacktrigon}$]

Lemmas \ref{noblackss}, \ref{noordi} and \ref{omega2} guarantee that
$\Lambda$ contains a $W$($\lambda_1, \lambda_2$) face where
$\{\lambda_1, \lambda_2\}$ = $\{\beta, \delta\}$ or $\{\alpha,
\gamma\}$ by the index equation $\sum _{vertices} I(v)$ +
$\sum_{faces} I(f)=2$. Also Lemmas \ref{noblackss}, \ref{noordi} and
\ref{omega1} guarantee that $\Lambda$ contains a $W$($\mu_1, \mu_2$)
face where $\{\mu_1, \mu_2\}$ = $\{\alpha, \delta\}$ or $\{\beta,
\gamma\}$. Hence a $W$($\lambda_1, \lambda_2$) face and a
$W$($\mu_1, \mu_2$) face are the only two non-isomorphic white disk
faces in $\Lambda$ by Theorem \ref{notriple}.

If $\{\lambda_1, \lambda_2\}$ = $\{\beta, \delta\}$ and $\{\mu_1,
\mu_2\}$ = $\{\alpha, \delta\}$, then no $\gamma$-edges are incident
to interior vertices, which is impossible since there are only four
types of interior vertices. See Figure 8. Similarly for the case
that $\{\lambda_1, \lambda_2\}$ = $\{\beta, \delta\}$ and $\{\mu_1,
\mu_2\}$ = $\{\beta, \gamma\}$, in which case no $\alpha$-edges are
incident to interior vertices, and the case that $\{\lambda_1,
\lambda_2\}$ = $\{\alpha, \gamma\}$ and $\{\mu_1, \mu_2\}$ =
$\{\alpha, \delta\}$, in which case no $\beta$-edges are incident to
interior vertices.

Then the only remaining case is that $\{\lambda_1, \lambda_2\}$ =
$\{\alpha, \gamma\}$ and $\{\mu_1, \mu_2\}$ = $\{\beta, \gamma\}$,
i.e. $\Lambda$ has a $W(\alpha, \gamma)$ face and a $W(\beta,
\gamma)$ face, in which case no $\delta$-edges are incident to
interior vertices, which is possible when an interior vertex is of
type ($i$) or ($ii$). These faces cannot be both bad by Lemma
\ref{order1} and cannot be both good by Lemma \ref{good}.
\begin{figure}[tbh]
\includegraphics{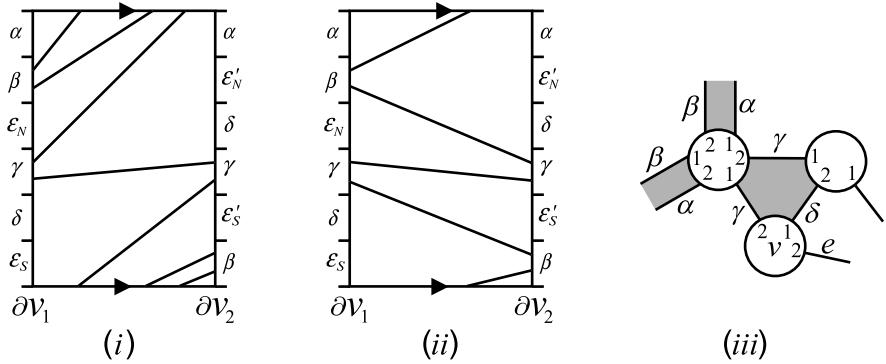}\caption{The intersection of the white corners of a bad $W(\beta,
\gamma)$ face and $H_W$, and the local configuration at an interior
vertex of type ($i$).}\label{f:jj}
\end{figure}

Assume that a $W(\alpha, \gamma)$ face is good and a $W(\beta,
\gamma)$ face is bad. Then every interior vertex is of type ($i$).
Consider the intersection of the white corners of a bad $W(\beta,
\gamma)$ face and $H_W$ where there are two possibilities as in
Figure 12 $(i), (ii)$, and the local configuration at an interior
vertex of type ($i$) as in Figure 12 $(iii)$. However, Figure 12
$(ii)$ is impossible since there is a $W(\alpha, \gamma)$ corner. In
Figure 12 ($iii$) $e$ must be $\gamma$ from Figure 12 $(i)$. Then the vertex $v$ has two
$\gamma$-edges at the same label, which is a contradiction to Lemma
\ref{same}.

Assume that a $W(\alpha, \gamma)$ face is bad and a $W(\beta,
\gamma)$ face is good. Then every interior vertex is of type ($ii$).
Consider the intersection of the white corners of a $W(\alpha,
\gamma)$ face and $H_W$ and the local configuration at an interior
vertex of type ($ii$). See Figure 13. There are two possibilities
for the intersection of the white corners of a $W(\alpha, \gamma)$
face and $H_W$ as in Figure 13 $(i), (ii)$. Figure 13 $(ii)$ is
impossible since there is a $W(\beta, \gamma)$ corner. In Figure 13
($iii$) suppose that $v_1$ is not an exceptional vertex. Then $e_1$ must
be an $\alpha$-edge and must be a boundary edge since the
$\delta$-edge adjacent to $e_1$ is a boundary edge. This implies
that $e_2$ is an interior edge and then a $\beta$-edge. Also, $e_3$
must be a $\gamma$-edge from Figure 13 $(i)$. This is a
contradiction to Lemma \ref{same}. If $v_1$ is an exceptional
vertex, consider the vertex $v_2$ in Figure 13
($iii$) and apply the same argument. Hence we
have finished the proof of Theorem \ref{noblacktrigon}.
\end{proof}

\begin{figure}[tbh]
\includegraphics{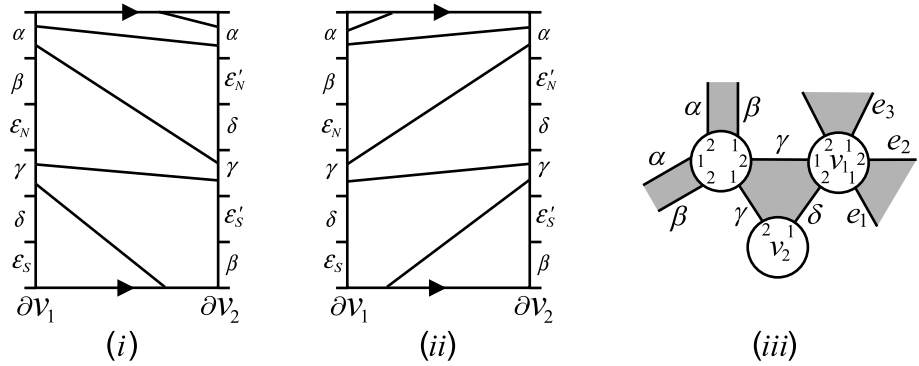}\caption{The intersection of the white corners of a $W(\alpha, \gamma)$
face and $H_W$, and the local configuration at an interior
vertex of type ($ii$).}\label{f:jj}
\end{figure}

\section{a $B(\alpha, \beta)$ bigon and a white trigon.}
In this section, we assume that $\Lambda$ contains a $B(\alpha,
\beta)$ bigon and does not contain a white bigon. Then the main goal of this section
is to prove the following theorem.
\begin{theorem} \label {nowhitetrigon}
If $\Lambda$ contains a $B(\alpha, \beta)$ bigon and does not contain a white bigon,
then $\Lambda$ does not contain a white trigon.
\end{theorem}

Suppose that $\Lambda$ contains a white trigon. Then Lemma
\ref{whit} says that a white trigon is either a $W(\alpha, \delta)$,
$W(\beta, \gamma)$, $W(\beta, \delta)$ or $W(\alpha, \gamma)$ trigon.

Note that Theorem \ref{notriple} still holds by using a white trigon
instead of a white binary disk face in the proof of Theorem
\ref{notriple}. Furthermore, by Lemma \ref{interior}, there is an
interior vertex in $\Lambda$.

First, assume that there is a $W(\alpha, \delta)$ trigon with two
$\delta$-edges and one $\alpha$-edge in $\Lambda$.
\begin{lemma} \label{twodeloneal1}
$\Lambda^*$ does not contain a boundary cycle with respect to
$\omega'$.
\end{lemma}
\begin{proof}
Consider the intersection of the white corners of a $W(\alpha,
\delta)$ trigon and $H_W$ where there are two cases as in Figure 14
$(i), (ii)$, and the intersection of the black corners of a
$B(\alpha, \beta)$ bigon and $H_B$ as in Figure 14 $(iii)$. Let $v$
be a boundary vertex dual to a boundary cycle with respect to
$\omega'$. There are two types of local configuration at $v$ as
shown in Figure 15 ($i$), ($ii$), where $e_1, e_2, e_3$ and $e_4$
are positive edges and $e_5$ and $e_6$ are negative edges.

\begin{figure}[tbh]
\includegraphics{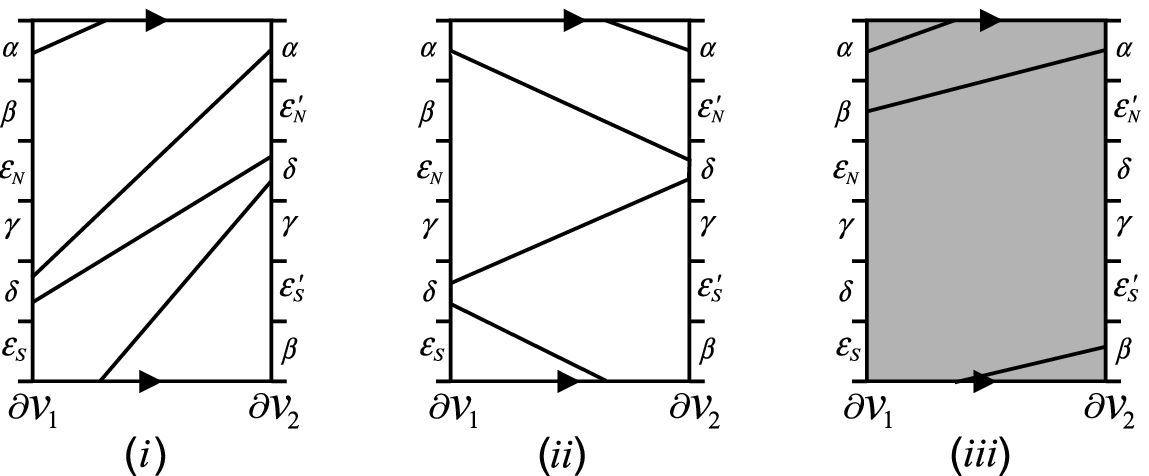}\caption{}\label{f:jj}
\end{figure}

\begin{figure}[tbh]
\includegraphics{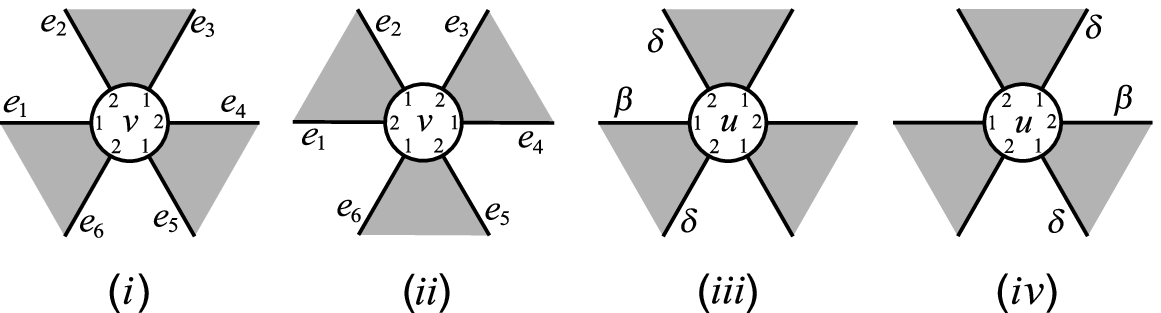}\caption{}\label{f:jj}
\end{figure}

Assume we have Figure 14 ($ii$). Then the only possible white corner
with one $\beta$-edge is a $W(\beta, \delta)$ corner or a $W(\delta,
\beta)$ corner, which makes a boundary vertex not to be a cycle with
respect to $\omega'$ since a $\beta$-edge and a $\delta$-edge have
the same orientation with respect to $\omega'$. Now assume Figure 14 ($i$).

1) The first type as in Figure 15 ($i$). Note that ($e_5, e_6)$ =
($\epsilon_N, \epsilon'_S$) from the intersection of the white
corners of a $W(\alpha, \delta)$ trigon and $H_W$ as shown in Figure
14 ($i$).

Assume there is a clockwise boundary cycle. There are only two cases for $e_1$
i.e. $e_1$ is either $\alpha$ or $\gamma$. First if $e_1$ =
$\alpha$, then $e_2$ must be $\delta$ from Figure 14 $(i)$ since
there is a $W(\epsilon_N, \epsilon'_S)$ corner. Then $e_3$ =
$\gamma$ and $e_4$ = $\beta$. Thus, we have a $B(\gamma, \delta)$
corner and a $B(\epsilon_N, \beta)$ corner. However this cannot
happen from the intersection of the black corners of a $B(\alpha,
\beta)$ bigon and $H_B$ as in Figure 14 ($iii$).

Second, if $e_1$ = $\gamma$, $e_2$ must be $\beta$ since there is a
$W(\epsilon_N, \epsilon'_S)$ corner. Then $e_3$ = $\alpha$ and $e_4$
= $\delta$. Hence we have two black corners, i.e. a $B(\gamma,
\epsilon'_S)$ corner and a $B(\epsilon_N, \delta)$ corner. We can
observe from Figure 14 $(iii)$ with a $B(\gamma, \epsilon'_S)$ corner
and a $B(\epsilon_N, \delta)$ corner inserted that the only positive
edge from which a black corner runs to a $\gamma$-edge is a
$\gamma$-edge. Thus if a $\gamma$-edge is an interior edge, then
there is a face all edges of which are $\gamma$-edges, which contradicts
Lemma 3.1 of \cite{GLu2}. Consequently, no
$\gamma$-edge is incident to interior vertices. For an interior
vertex $u$, a $\delta$-edge must be incident to $u$, otherwise there
are three $\alpha$ or $\beta$-edges incident to $u$, which is a
contradiction to Lemma \ref{same}. However at whatever label a
$\delta$-edge is incident to $u$, there is another $\delta$-edge
incident to $u$ at the same label from Figure 14 $(i)$ with a
$W(\epsilon_N, \epsilon'_S)$ corner inserted and Figure 14 $(iii)$ with a $B(\gamma, \epsilon'_S)$ corner
and a $B(\epsilon_N, \delta)$ corner inserted, which is a
contradiction to Lemma \ref{same}. See Figure 15 ($iii$), ($iv$).

Assume there is a anticlockwise boundary cycle. If $e_1$ = $\beta$,
then $e_2$ = $\gamma$, $e_3$ = $\delta$ and $e_4$ = $\alpha$ since
there is a $W(\epsilon_N, \epsilon'_S)$ corner. Thus, we have a
$B(\delta, \gamma)$ corner and a $B(\beta, \epsilon'_S)$ corner.
However these two black corners cannot happen at the same time from
the intersection of the black corners of a $B(\alpha, \beta)$ bigon
and $H_B$ as in Figure 14 $(iii)$. If $e_1$ = $\delta$, then $e_2$ =
$\alpha$, $e_3$ = $\beta$ and $e_4$ = $\gamma$. Hence we have two
black corners, i.e. a $B(\epsilon_N, \gamma)$ corner and a
$B(\delta, \epsilon'_S)$ corner at $v$. We can observe that the only
positive edge class to which a black corner runs from a
$\gamma$-edge is a $\gamma$-edge. The argument is then completely
the same as that of the case of a clockwise boundary cycle.

2) The second type as in Figure 15 ($ii$). We use a similar argument
as in the proof of Lemma \ref{omega2}. If we let $u_x$ be a boundary
vertex dual to a boundary cycle, then the label $x$ appears either three
times around $v_1$ at the ends of $\alpha$-, $\gamma$- and
$\epsilon$-edges and three times around $v_2$ at the ends of
$\beta$-, $\delta$- and $\epsilon'$-edges, or three times around
$v_1$ at the ends of $\beta$-, $\delta$- and $\epsilon$-edges and
three times around $v_2$ at the ends of $\alpha$-, $\gamma$- and
$\epsilon'$-edges. By the Claim in Lemma \ref{omega2}, $(e_5, e_6)$
is either ($\epsilon_N, \epsilon'_S$) or ($\epsilon_S,
\epsilon'_N$). However, one of $e_1$ and $e_4$ must be either
$\alpha$ or $\gamma$. Thus by Figure 14 $(i)$, $(e_5, e_6)$ =
($\epsilon_N, \epsilon'_S$). Figure 14 ($iii$) with a $B(\epsilon_N,
\epsilon'_S)$ corner inserted shows that if one edge of a black
corner is a $\delta$-edge, then the other of the corner must be a
$\beta$-edge. Therefore there is no boundary cycle with respect to
$\omega'$. This completes the proof.
\end{proof}

\begin{figure}[tbh]
\includegraphics{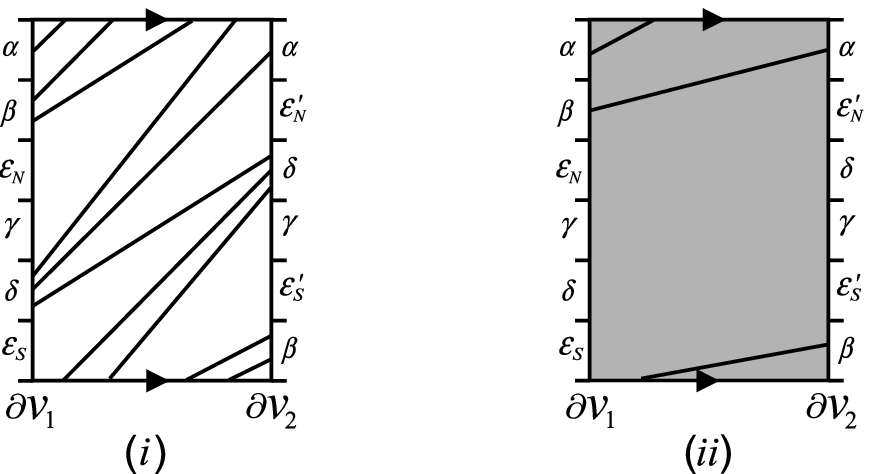}\caption{}\label{f:jj}
\end{figure}

\begin{figure}[tbh]
\includegraphics{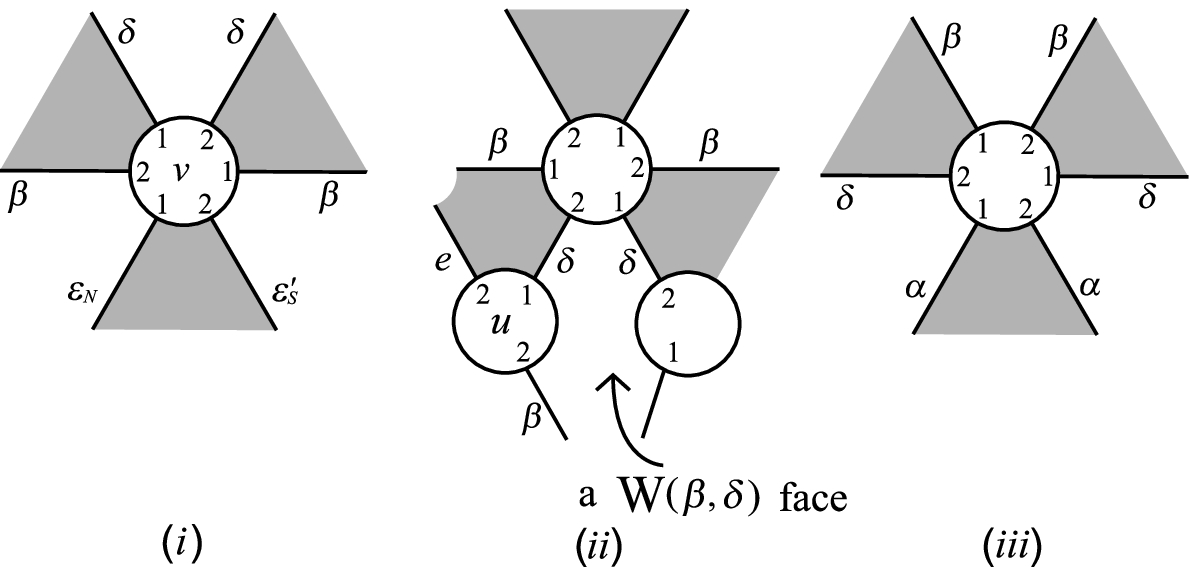}\caption{}\label{f:jj}
\end{figure}

\begin{lemma} \label{twodeloneal2}
There is no source or sink dual to a white face in $\Lambda^*$ with
respect to $\omega'$.
\end{lemma}
\begin{proof}
Suppose there is a source or sink dual to a white face with respect
to $\omega'$, i.e. there is a $W(\lambda, \mu)$ face where
$\{\lambda, \mu\}$ = $\{\alpha, \gamma\}$ or $\{\beta, \delta\}$.
Since $\Lambda$ contains a $W(\alpha, \delta)$ trigon, which is
good, Lemma \ref{good} implies that a $W(\lambda, \mu)$ face is bad.
If $\{\lambda, \mu\}$ = $\{\alpha, \gamma\}$ i.e. $\Lambda$ contains
a bad $W(\alpha, \gamma)$ face, then we have a $W(\alpha, \alpha)$
corner, a $W(\gamma, \gamma)$ corner and a $W(\delta, \delta)$
corner. This is impossible by Lemma \ref{order1}. Therefore we
assume that $\Lambda$ has a bad $W(\beta, \delta)$ face. Note that
Theorem \ref{notriple} implies that every white face in $\Lambda$ is
either a $W(\alpha, \delta)$ trigon or a bad $W(\beta, \delta)$
face. Consider the intersection of the white corners of a $W(\alpha,
\delta)$ trigon, the white corners of a bad $W(\beta, \delta)$ face
and $H_W$ where there is only one case, and the intersection of the
black corners of a $B(\alpha, \beta)$ bigon and $H_B$. See Figure
16. Since $\Lambda$ has a $W(\alpha, \delta)$ trigon and a bad
$W(\beta, \delta)$ face, there are at least two vertices which
have a $W(\delta, \delta)$ corner. Let $v$ be such a non-exceptional
vertex.

\begin{claim}
There is only one local configuration at $v$ as shown in Figure $17$
$(i)$.
\end{claim}

\begin{proof}
If at least five positive edges are incident to $v$, then there are
at least two white corners at $v$ which are corners of white disk
faces in $\Lambda$. One is a $W(\delta, \delta)$ corner. The other
corner must be a $W(\beta, \beta)$ corner because every white face
in $\Lambda$ is either a $W(\alpha, \delta)$ trigon or a bad
$W(\beta, \delta)$ face and no vertex has three edges incident to it
in the same edge class. The one remaining white corner at $v$ must
be a $W(\alpha/\gamma/\epsilon, \alpha/\gamma/\epsilon')$ corner,
which is impossible by Figure 16 ($i$). If four edges incident to
$v$ are positive, the local configuration at $v$ must be as in
Figure 17 ($i$) since a $W(\epsilon, \epsilon')$ corner cannot occur
by Figure 16 ($i$).
\end{proof}

Suppose a $W(\delta, \delta)$ corner of a bad $W(\beta, \delta)$
face doesn't appear at an exceptional vertex. The bad $W(\beta,
\delta)$ face contains a vertex $u$ as shown in Figure 17 ($ii$).
Since there is only one local configuration at the vertex containing
a $W(\delta, \delta)$ corner, $e$ is positive. However, $e$ cannot
be $\alpha, \gamma$ and $\delta$ from Figure 16 $(ii)$ because there
is a $B(\epsilon_N, \epsilon'_S)$ corner by the Claim. Also, $e$
cannot be $\beta$ by Lemma \ref{same}.

Finally, every $W(\delta, \delta)$ corner of a $W(\beta, \delta)$
face must appear at an exceptional vertex. Therefore there is only
one $W(\beta, \delta)$ face $X$, say. Thus there is only one
possible sink or source dual to a white face. Note that an
exceptional vertex is not dual to a cycle in $\Lambda^*$ because it
contains a $W(\delta, \delta)$ corner. On the other hand, every white face in $\Lambda$ is
either a $W(\alpha, \delta)$ trigon or a bad $W(\beta, \delta)$
face. Thus there is a
unique type of interior vertices as in Figure 17 ($iii$) since a
$B(\epsilon_N, \epsilon'_S)$ corner coming from the Claim and Figure
16 ($ii$) exclude a $B(\delta, \delta)$ corner. Since there are
only two non-isomorphic black disk faces in $\Lambda$, one black
disk face has an $\alpha$-, $\beta$- and $\delta$-edges coming from
the local picture at an interior vertex as shown in Figure 17 ($iii$).
The other black disk face is a $B(\alpha, \beta)$ bigon. Hence there
is no sink or source dual to a black disk face. Since there is no
boundary cycle by Lemma \ref{twodeloneal1}, the only possible
positive index takes place at $X$ or the outside face. However, the
index of an interior vertex is $-1$. This is a contradiction to the
equation $\sum _{vertices} I(v)$ + $\sum_{faces} I(f)=2$.
\end{proof}

\begin{lemma} \label{twodeloneal3}
There is no source or sink dual to a black face in $\Lambda^*$ with
respect to $\omega'$.
\end{lemma}

\begin{figure}[tbh]
\includegraphics{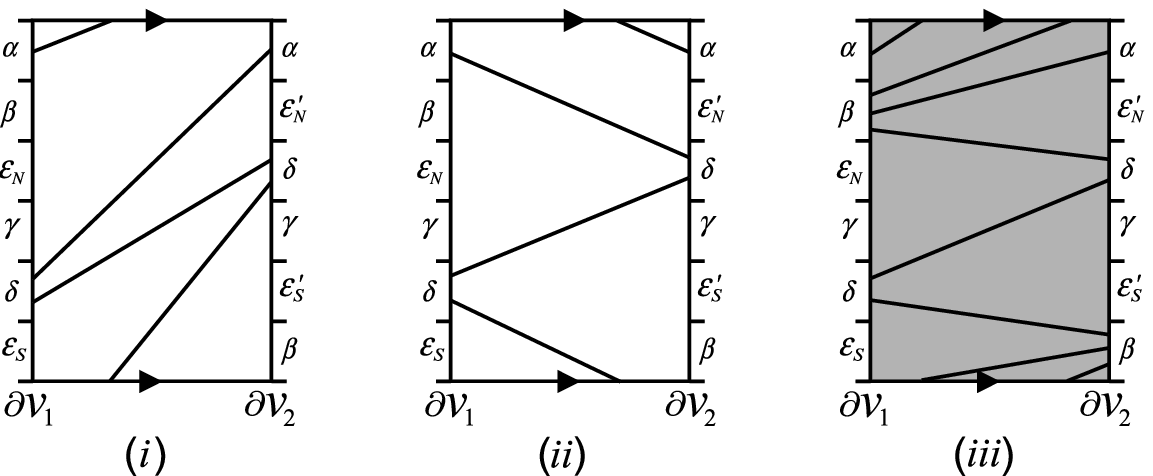}\caption{}\label{f:jj}
\end{figure}

\begin{figure}[tbh]
\includegraphics{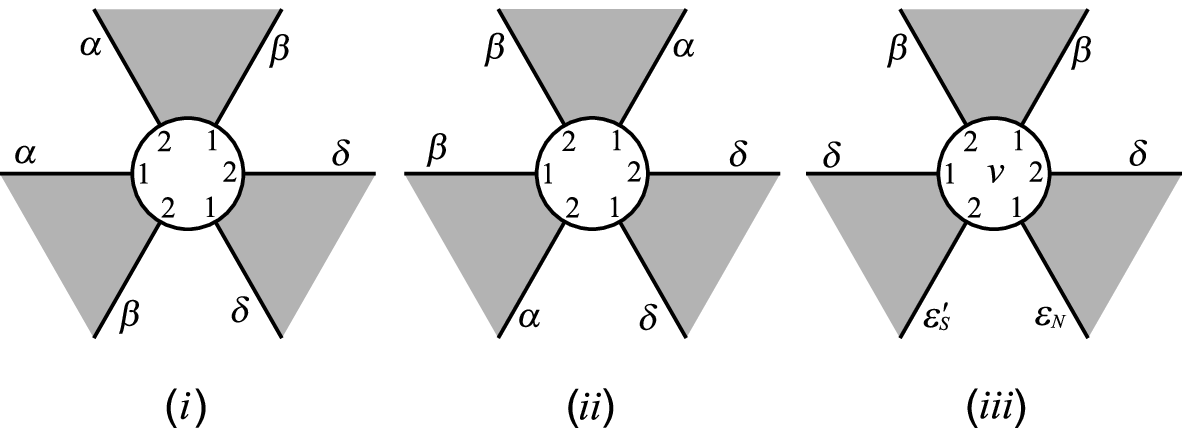}\caption{}\label{f:jj}
\end{figure}

\begin{proof}
Suppose there is a source or sink dual to a black face with respect
to $\omega'$, i.e. there is a $B(\lambda, \mu)$ face where
$\{\lambda, \mu\}$ = $\{\alpha, \gamma\}$ or $\{\beta, \delta\}$.
Then a $B(\lambda, \mu)$ face must be bad by Lemma \ref{good}
because of the existence of a $B(\alpha, \beta)$ bigon.

Suppose we have a bad $B(\beta, \delta)$ face. Consider the
intersection of the black corners of a $B(\alpha, \beta)$ bigon, the
black corners of a bad $B(\beta, \delta)$ face and $H_B$ as
shown in Figure 18 $(iii)$, and the intersection of the white
corners of a $W(\alpha, \delta)$ trigon and $H_W$ where there are two
cases as shown in Figure 18 $(i), (ii)$. By Theorem \ref{notriple},
every black face is either a $B(\alpha, \beta)$ bigon or a bad
$B(\beta, \delta)$ face. Also, since three edges in the same edge
class cannot be incident to a vertex, there are only two types of
local configuration at interior vertices. See Figure 19 ($i$)
($ii$). Note that Figure 18 ($i$) excludes the first type of
interior vertices because a $W(\alpha, \alpha)$ corner and a
$W(\delta, \beta)$ corner cannot exist simultaneously in Figure 18
($i$), and Figure 18 ($ii$) rules out the second type of interior
vertices since no $W(\beta, \beta)$ corner survives in Figure 18
($ii$).

Let $v$ be a vertex containing a $B(\beta, \beta)$ corner of a bad
$B(\beta, \delta)$ face. Suppose $v$ is a non-exceptional vertex.
Now apply a similar argument as in the proof of the Claim in Lemma
\ref{twodeloneal2} by using Figure 18 ($iii$). Then the local configuration at $v$ looks like
Figure 19 ($iii$). However, a $W(\epsilon_N, \epsilon'_S)$ corner is
impossible because in the case of Figure 18 ($i$) the second type of
interior vertices contains a $W(\beta, \beta)$ corner, and in the
case of Figure 18 ($ii$) a $W(\epsilon_N, \epsilon'_S)$ corner
cannot exist.

Consequently, $v$ must be an exceptional vertex and then there is
only one bad $B(\beta, \delta)$ face $Y$, say, containing only one
$B(\beta, \beta)$ corner. Note that since $v$ contains a $B(\beta,
\beta)$ corner, the exceptional vertex $v$ is not dual to a cycle in
$\Lambda^*$. We already showed that there are no $W(\beta, \delta)$
faces and no $W(\alpha, \gamma)$ faces in $\Lambda$. Furthermore, no
$B(\alpha, \gamma)$ faces exist since we have a bad $B(\beta,
\delta)$ face. Therefore, there are only two possible sinks or
sources in $\Lambda^*$. One is dual to $Y$ and the other the outside
face.

Since there is no boundary cycle in $\Lambda^*$ by Lemma \ref{twodeloneal1},
in order to get a contradiction to the equation $\sum _{vertices}
I(v)$ + $\sum_{faces} I(f)=2$, it suffices to show that either there
is a vertex of index $-1$ or the outside face is not dual to a sink
or source. Let $u$ be a vertex containing a $B(\beta, \delta)$
corner of $Y$. Then $u$ is not an exceptional vertex. Also $u$ is a
boundary vertex since the local picture at interior vertices looks like Figure 19
($i$) ($ii$). We can observe that two negative edges are
not incident to $u$ since neither a $W(\epsilon, \epsilon')$ corner
nor a $B(\epsilon, \epsilon')$ corner exists by Figure 18 with white
corners of an interior vertex added. If at most one negative edge is
incident to $u$, it is easy to see using Figure 18 that either $u$
has index $-1$ or the outside face is not dual to a sink or a source.

Suppose we have a bad $B(\alpha, \gamma)$ face in $\Lambda$. Apply
a similar argument as in the case of a $B(\beta, \delta)$ face.
First use a $B(\alpha, \alpha)$ corner instead of a $B(\beta,
\beta)$ corner to show that there is only one $B(\alpha, \gamma)$
face in $\Lambda$. However unlike the case of a $B(\beta, \delta)$
face, there is only one type of local configuration at interior
vertices and furthermore the index of an interior vertex is $-1$. Use
the equation $\sum _{vertices} I(v)$ + $\sum_{faces} I(f)=2$ to get
a contradiction.
\end{proof}

By Lemmas \ref{twodeloneal1}, \ref{twodeloneal2}, and
\ref{twodeloneal3} and the index equation $\sum _{vertices} I(v)$ +
$\sum_{faces} I(f)=2$, we have proved Theorem \ref{nowhitetrigon}
for the case that $\Lambda$ contains a $W(\alpha, \delta)$ trigon
with two $\delta$-edges and one $\alpha$-edge.

Next, suppose $\Lambda$ contains a $W(\alpha, \delta)$ trigon with
two $\alpha$-edges and one $\delta$-edge.

\begin{lemma} \label{twoalonedel1}
$\Lambda^*$ does not contain a boundary cycle with respect to
$\omega'$.
\end{lemma}
\begin{proof}
The proof of this lemma is exactly analogous as that of Lemma
\ref{twodeloneal1}.
\end{proof}
\begin{lemma} \label{twoalonedel2}
There is no source or sink dual to a white face in $\Lambda^*$ with
respect to $\omega'$.
\end{lemma}

\begin{figure}[tbh]
\includegraphics{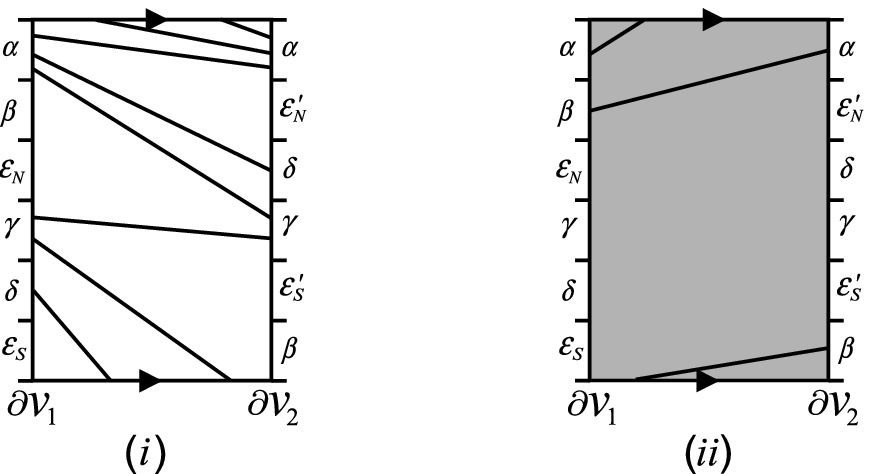}\caption{}\label{f:jj}
\end{figure}

\begin{figure}[tbh]
\includegraphics{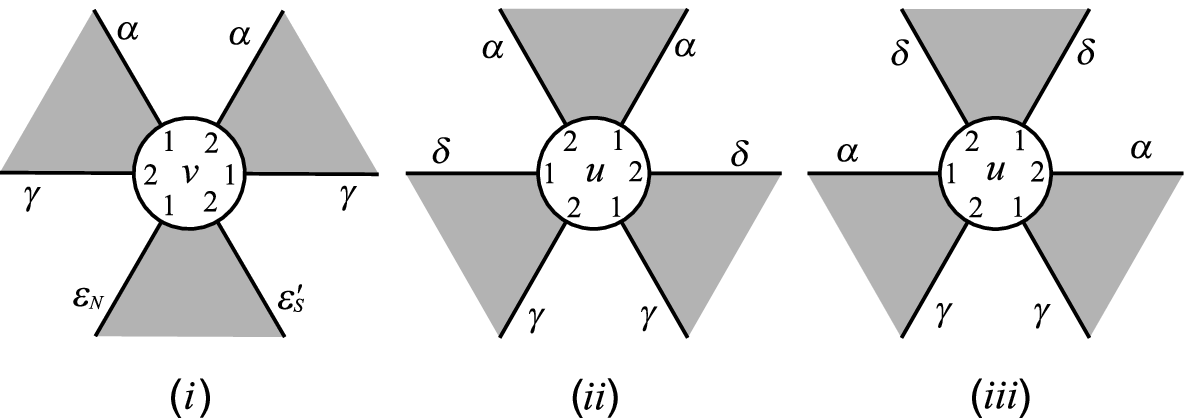}\caption{}\label{f:jj}
\end{figure}

\begin{proof}
Suppose there is a source or sink dual to a white face with respect
to $\omega'$, i.e. there is a $W(\lambda, \mu)$ face where
$\{\lambda, \mu\}$ = $\{\alpha, \gamma\}$ or $\{\beta, \delta\}$.
Then a $W(\lambda, \mu)$ face must be a bad $W(\alpha, \gamma)$ face
by Lemma \ref{good} and Lemma \ref{order1}. Consider the
intersection of the white corners of a $W(\alpha, \delta)$ trigon,
the white corners of a bad $W(\alpha, \gamma)$ face and $
H_W$ where there is only one case, and the intersection of the black
corners of a $B(\alpha, \beta)$ bigon and $H_B$. See Figure
20. Since $\Lambda$ contains a $W(\alpha, \delta)$ trigon and a bad
$W(\alpha, \gamma)$ face, there are at least two vertices which have
a $W(\alpha, \alpha)$ corner. Let $v$ be such a non-exceptional
vertex.

\begin{claim}
There is only one local configuration at $v$ as shown in Figure $21$
$(i)$.
\end{claim}
\begin{proof}
Follow the proof of the Claim in Lemma \ref{twodeloneal2} with
$\beta, \delta$ replaced by $\gamma, \alpha$ respectively.
\end{proof}

Let $u$ be an interior vertex in $\Lambda$. Since every
white face in $\Lambda$ is either a $W(\alpha, \delta)$ trigon or a
bad $W(\alpha, \gamma)$ face by Theorem \ref{notriple}, there are only
two possible types of the local configuration at $u$ by Lemma
\ref{same}. See Figure 21 ($ii$), ($iii$). However either type is
also impossible since a $B(\delta, \gamma)$ corner in Figure 21
($ii$) and a $B(\delta, \delta)$ corner in Figure 21 ($iii$) cannot
happen by the existence of a $B(\epsilon_N, \epsilon'_S)$ corner
coming from the local configuration at $v$ in the Claim. This
completes the proof.
\end{proof}

\begin{lemma} \label{twoalonedel3}
There is no source or sink dual to a black face in $\Lambda^*$ with
respect to $\omega'$.
\end{lemma}
\begin{proof}
Suppose there is a source or sink dual to a black face with respect
to $\omega'$, i.e. there is a $B(\lambda, \mu)$ face where
$\{\lambda, \mu\}$ = $\{\alpha, \gamma\}$ or $\{\beta, \delta\}$.
Then a $B(\lambda, \mu)$ face must be bad by Lemma \ref{good} since
$\Lambda$ contains a $B(\alpha, \beta)$ bigon.

\begin{figure}[tbh]
\includegraphics{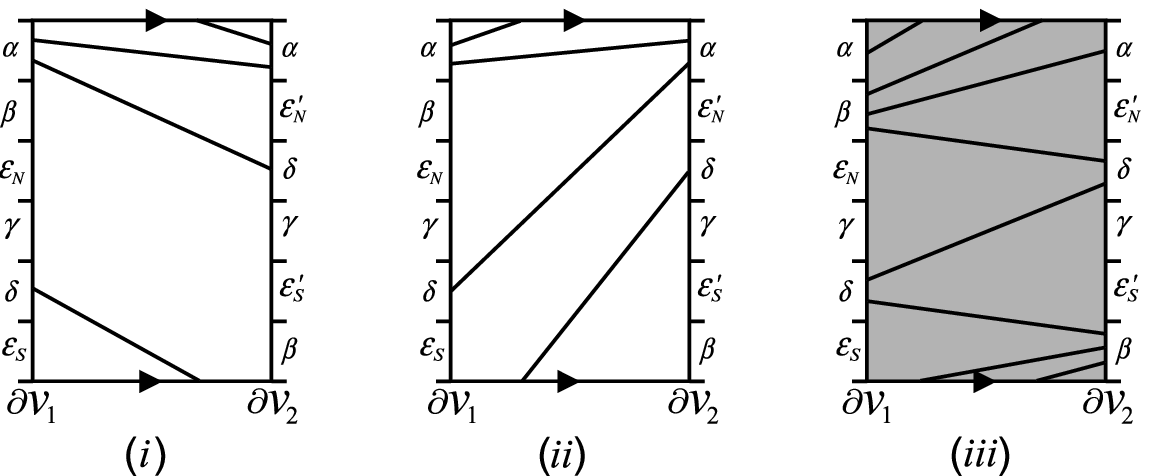}\caption{}\label{f:jj}
\end{figure}

\begin{figure}[tbh]
\includegraphics{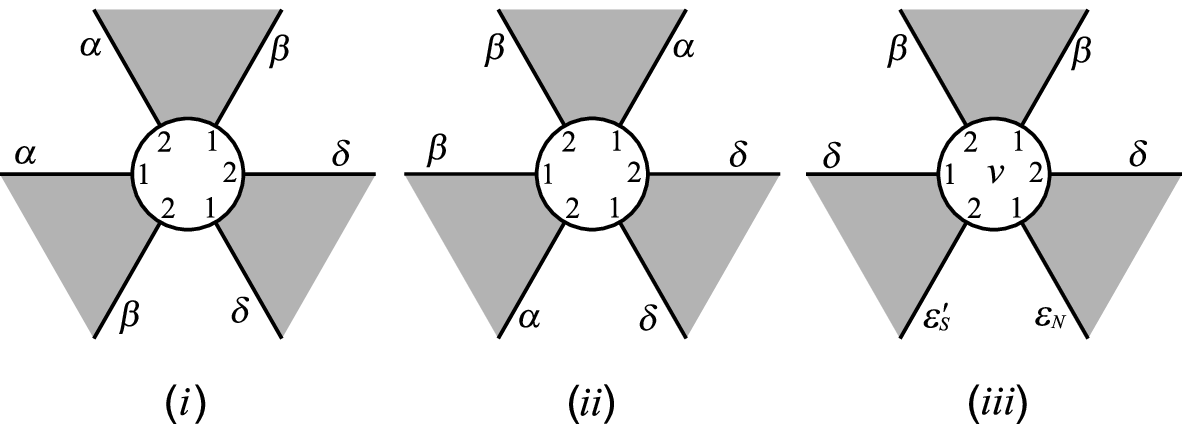}\caption{}\label{f:jj}
\end{figure}

Suppose we have a bad $B(\beta, \delta)$ face. Consider the
intersection of the black corners of a $B(\alpha, \beta)$ bigon, the
black corners of a bad $B(\beta, \delta)$ face and $H_B$ as
in Figure 22 $(iii)$, and the intersection of the white corners of a
$W(\alpha, \delta)$ trigon and $H_W$ where there are two
cases as in Figure 22 $(i), (ii)$. By Theorem \ref{notriple}, every
black face is either a $B(\alpha, \beta)$ bigon or a bad $B(\beta,
\delta)$ face. Also, since three edges in the same edge class cannot
be incident to a vertex, there are only two types of the local
configuration at interior vertices. See Figure 23 ($i$), ($ii$).
Note that the white corners of an interior vertex rule out Figure 22
($ii$).
\begin{claim}
A bad $B(\beta, \delta)$ face whose $B(\beta, \beta)$ corner appears
at a non-exceptional vertex contributes $0$ to the sum $\sum
_{vertices} I(v)$ + $\sum_{faces} I(f)$.
\end{claim}
\begin{proof}
Let $v$ be a vertex containing a $B(\beta, \beta)$ corner of a bad
$B(\beta, \delta)$ face. A similar argument as in the proof of the
Claim in Lemma \ref{twodeloneal2} shows that the local configuration at
$v$ looks like Figure 23 ($iii$). Then $v$ has index $-1$. This
completes the proof.
\end{proof}

Suppose no $B(\beta, \beta)$ corners appear at an exceptional
vertex. Lemmas \ref{twoalonedel1}, \ref{twoalonedel2} and the Claim
imply that the only possible positive index comes from the cycle
dual to an exceptional vertex or the sink or source dual to the
outside face. However both cannot happen simultaneously so we get a
contradiction to the equation $\sum _{vertices} I(v)$ +
$\sum_{faces} I(f)=2$. Hence there is a $B(\beta, \beta)$ corner of
a $B(\beta, \delta)$ face which must appear at an exceptional
vertex. Then the rest of the proof is pretty similar to the fourth and the fifth
paragraphes of the proof of Lemma \ref{twodeloneal3}.

For the case that there is a bad $B(\alpha, \gamma)$ face in
$\Lambda$, apply a similar argument as in the case of a
$B(\beta, \delta)$ face.
\end{proof}

By Lemmas \ref{twoalonedel1}, \ref{twoalonedel2}, and
\ref{twoalonedel3} and the index equation $\sum _{vertices} I(v)$ +
$\sum_{faces} I(f)=2$, we have proved Theorem \ref{nowhitetrigon}
for the case that $\Lambda$ contains a $W(\alpha, \delta)$ with two
$\alpha$-edges and one $\delta$-edge.

For the case of a $W(\beta, \gamma)$ trigon, we can apply the same
argument as in the case of a $W(\alpha, \delta)$ trigon by
interchanging $\alpha$ and $\beta$, and by interchanging $\gamma$
and $\delta$.

Now suppose $\Lambda$ contains a $W(\beta, \delta)$ trigon with two
$\delta$-edges and one $\beta$-edge. Here, instead of using the
orientation $\omega'$, we will use the orientation $\omega$.

\begin{figure}[tbh]
\includegraphics{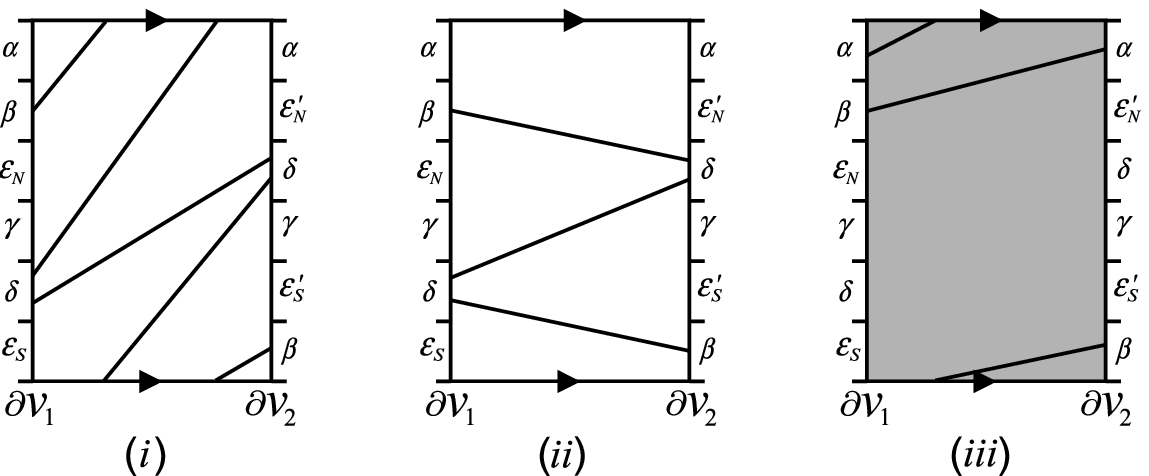}\caption{}\label{f:jj}
\end{figure}

\begin{figure}[tbh]
\includegraphics{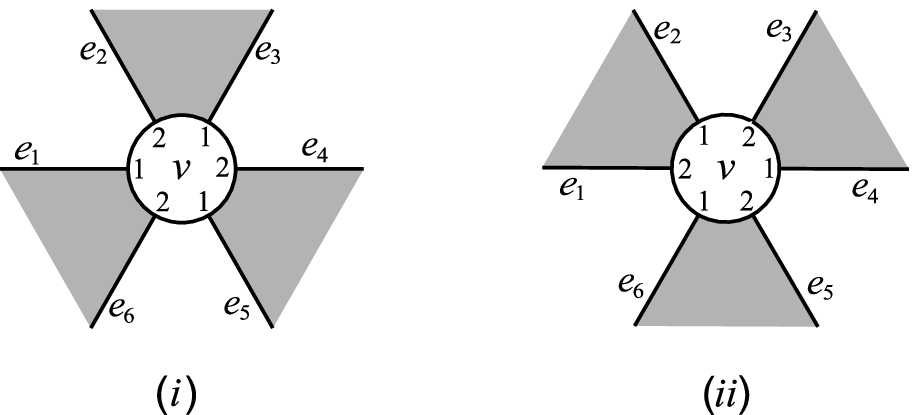}\caption{}\label{f:jj}
\end{figure}

\begin{lemma} \label{twodelonebe1}
$\Lambda^*$ does not contain a boundary cycle with respect to
$\omega$.
\end{lemma}
\begin{proof}
Consider the intersection of the white corners of a $W(\beta,
\delta)$ trigon and $H_W$, where there are two cases as in
Figure 24 $(i), (ii)$, and the intersection of the black corners of
a $B(\alpha, \beta)$ bigon and $H_B$ as in Figure 24
$(iii)$. Let $v$ be a boundary vertex dual to a boundary cycle with
respect to $\omega$. There are two types of the local configuration at $v$ as
shown in Figure 25 where $e_1, e_2, e_3$ and $e_4$ are positive
edges, and $e_5$ and $e_6$ are negative edges.

First we assume Figure 24 ($i$); then if one edge of a white corner
is an $\alpha$-edge, then the other edge of the corner must be a
$\delta$-edge. Since an $\alpha$-edge and a $\delta$-edge have the
same orientation with respect to $\omega$, there is no boundary cycle, as desired.

Second we assume Figure 24 ($ii$); then

1) The first type as in Figure 25 ($i$);  Note that ($e_5, e_6)$ =
($\epsilon_S, \epsilon'_N$) from the intersection of the white
corners of a $W(\beta, \delta)$ trigon and $H_W$ as in
Figure 24 ($ii$). That is, there is a $W(\epsilon_S, \epsilon'_N)$
corner. Then Figure 24 ($ii$) with this corner inserted shows that
if one edge of a white corner in $\Lambda$ is an $\alpha$-edge, then
the other edge must be a $\delta$-edge. Therefore, there is no
boundary cycle.

2) The second type as in Figure 25 ($ii$); Assume there is a
clockwise boundary cycle. If $e_2$ = $\gamma$ and $e_4$ = $\beta$,
then $e_1$ = $\alpha$ and $e_3$ = $\delta$. In other words, we have
a $B(\gamma, \alpha)$ corner and a $B(\beta, \delta)$ corner, which
is impossible by Figure 24 ($iii$). If $e_2$ = $\beta$ and $e_4$ =
$\gamma$, then $e_5$ must be $\delta$ by Figure 24 ($ii$), which is
absurd because $e_5$ is a negative edge.

Assume there is a anticlockwise boundary cycle. If $e_1$ = $\gamma$,
$e_6$ must be $\delta$ by Figure 24 ($ii$), which is absurd. If
$e_1$ = $\beta$ and $e_3 = \gamma$, then $e_2$ = $\delta$ and $e_4$
= $\alpha$. In other words, we have a $B(\alpha, \gamma)$ corner and
a $B(\delta, \beta)$ corner, which is impossible by Figure 24
($iii$).
\end{proof}

\begin{lemma} \label{twodelonebe2}
There is no source or sink dual to a white face in $\Lambda^*$ with
respect to $\omega$.
\end{lemma}
\begin{proof}
Apply exactly the same argument as that of Lemma
\ref{twodeloneal2} by interchanging $\alpha$ and $\beta$.
\end{proof}

\begin{lemma} \label{twodelonebe3}
There is no source or sink dual to a black face in $\Lambda^*$ with
respect to $\omega$.
\end{lemma}
\begin{proof}
The proof is exactly the same as that of Lemma \ref{twodeloneal3}
with $\alpha$ and $\beta$ interchanged.
\end{proof}

In conclusion, by Lemmas \ref{twodelonebe1}, \ref{twodelonebe2}, and
\ref{twodelonebe3} and the index equation $\sum _{vertices} I(v)$ +
$\sum_{faces} I(f)=2$, we have proved Theorem \ref{nowhitetrigon}
for the case that $\Lambda$ contains a $W(\beta, \delta)$ trigon
with two $\delta$-edges and one $\beta$-edge.

Now suppose $\Lambda$ contains a $W(\beta, \delta)$ with two
$\beta$-edges and one $\delta$-edge.

\begin{lemma} \label{twobeonedel1}
$\Lambda^*$ does not contain a boundary cycle with respect to
$\omega$.
\end{lemma}
\begin{proof}
The proof is completely similar as that of Lemma \ref{twodelonebe1}.
\end{proof}

\begin{lemma} \label{twobeonedel2}
There is no source or sink dual to a white face in $\Lambda^*$ with
respect to $\omega$.
\end{lemma}
\begin{proof}
The proof is exactly the same as that of Lemma \ref{twoalonedel2}
with $\alpha$ and $\beta$ interchanged.
\end{proof}

\begin{lemma} \label{twobeonedel3}
There is no source or sink dual to a black face in $\Lambda^*$ with
respect to $\omega$.
\end{lemma}
\begin{proof}
The proof is exactly the same as that of Lemma \ref{twoalonedel3}
with $\alpha$ and $\beta$ interchanged.
\end{proof}

In conclusion, by Lemmas \ref{twobeonedel1}, \ref{twobeonedel2}, and
\ref{twobeonedel3} and the equation $\sum _{vertices} I(v)$ +
$\sum_{faces} I(f)=2$, we have proved Theorem \ref{nowhitetrigon}
for the case that $\Lambda$ contains a $W(\beta, \delta)$ trigon
with two $\beta$-edges and one $\delta$-edge.

For the case of a $W(\alpha, \gamma)$ trigon, we can apply the same
argument by interchanging $\alpha$ and $\beta$, and by exchanging
$\gamma$ and $\delta$.

Finally, this completes the proof of Theorem \ref{nowhitetrigon}.

Theorem \ref{whitebigon} says that $G_1$ contains a $B(\alpha,
\beta)$ face and a $W(\lambda, \mu)$ bigon, where $\{\lambda, \mu
\}$ = $\{\alpha, \delta \}$, $\{\alpha, \gamma \}$, $\{\beta, \delta
\}$ or $\{\beta, \gamma \}$. We may assume without loss of
generality that $G_1$ contains a $B(\alpha, \beta)$ face and a
$W(\alpha, \delta)$ bigon.

The next proposition will be used to prove Theorem \ref{23}.
\begin{proposition} \label{trigon}
$G_1$ contains either a black bigon with two positive edges or a
trigon.
\end{proposition}

$Remark$: There are two types of trigons in $G_1$; one is a trigon
with all positive edges and the other is a trigon with one positive
edge and two negative edges.

To prove Proposition \ref{trigon}, we need two lemmas.

\begin{lemma} \label{negativebigon}
There is no bigon in $G_1$ whose edges are negative with a $(\epsilon_N,
\epsilon'_S)$ corner and a $(\epsilon_S, \epsilon'_N)$ corner as
shown in Figure $26$ $(i)$.
\end{lemma}

\begin{figure}[tbh]
\includegraphics{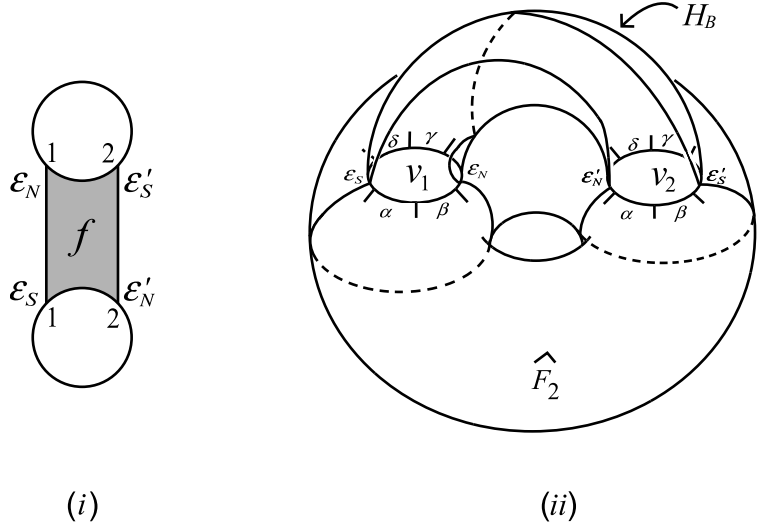}\caption{A black bigon with negative edges and with a $(\epsilon_N,
\epsilon'_S)$ corner and a $(\epsilon_S, \epsilon'_N)$ corner.}\label{f:jj}
\end{figure}

\begin{proof} Let $f$ be such a bigon as in the hypothesis. We assume
without loss of generality that $f$ is black. Recall that
$\widehat{F}_2$ is an incompressible torus containing the graph
$G_2$ and $H_B$ is the intersection of $V_2$ and a black
side of $M$. $\partial f$ appears on $\widehat{F}_2 \cup
H_B$ as in Figure 26 ($ii$).

The two endpoints of the negative edges $\epsilon$ and $\epsilon'$
of $f$ separate $\partial v_1$ and $\partial v_2$ into two arcs
respectively. Let $a$ (resp. $b$) be the arc on $\partial v_1$
(resp. $\partial v_2$) which contains the endpoints of $\alpha$-,
$\beta$-edges at $v_1$ (resp. $\gamma$-, $\delta$-edges at $v_2$). Then
$\epsilon$ $\cup$ $a$ and $\epsilon'$ $\cup$ $b$ are circles which
cobound an annulus on $\widehat{F}_2$, denoted by $A$. Also the two
corners of $f$ separate $\overline{H_B - (v_1 \cup v_2)}$
into two disks; let $D$ be one which contains $a$ and $b$. Then $A
\cup D$ is an once punctured Klein bottle. Therefore $A \cup D \cup
f$ is a Klein bottle, which is a contradiction.
\end{proof}

\begin{lemma} \label{nowhitecorner}
No white $(\epsilon_N, \epsilon'_N)$ or $(\epsilon_S, \epsilon'_S)$
corners exist in $G_1$.
\end{lemma}
\begin{proof}
Since there is a $W(\alpha, \delta)$ bigon, we have a $W(\alpha,
\delta)$ corner and a $W(\delta, \alpha)$ corner. Then it follows
from Lemma \ref{order} that a white $(\epsilon_N, \epsilon'_N)$ or
$(\epsilon_S, \epsilon'_S)$ corner cannot exist.
\end{proof}
\begin{proof}[Proof of Proposition $\ref{trigon}$]
Suppose for contradiction that $G_1$ does not contain a black bigon with
positive edges and does not contain a trigon. Observe that Lemmas
\ref{negativebigon}, \ref{nowhitecorner} imply that
$G_1$ doesn't contain a white bigon with two negative edges. Let
$\overline{G}$ be a connected component of the reduced graph
$\overline{G}_1$ of $G_1$. Note that $\overline{G}$ is a graph on
the sphere $\widehat{F}_1$.
\begin{claim}
Every vertex in $\overline{G}$ has valency at least $4$.
\end{claim}
\begin{proof}
Since $G_1$ does not contain a black bigon with positive edges and does not contain a white
bigon with negative edges, it is easy to see that no vertex in
$\overline{G}$ has valency 1 or 2. Suppose a vertex $v$ in
$\overline{G}$ has valency 3. Then there are three types of the local
configuration at $v$ in $G_1$. See Figure 27.

\begin{figure}[tbh]
\includegraphics{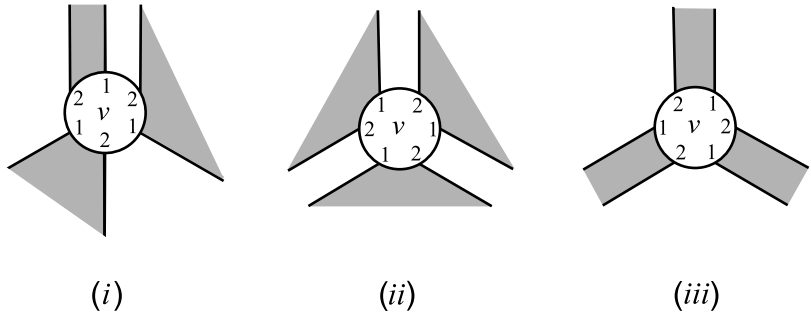}\caption{A vertex of valency 3 in $\overline{G}$.}\label{f:jj}
\end{figure}
For the first type, since there is no black bigon with positive
edges, the edges of the black bigon in Figure 27 ($i$) are negative.
Thus the adjacent white bigon has two negative edges, which is a
contradiction.

For the second type, all three white bigons in Figure 27 ($ii$) are
$W(\alpha, \delta)$ bigons since there is no white bigon with
negative edges in $G_1$. Therefore there are two edges in the same
edge class at the same label at $v$, which is a contradiction to
Lemma \ref{same}.

If we have the third type as in Figure 27 $(iii)$, all black bigons
have negative edges. By Lemma \ref{negativebigon}, all the corners of
these black bigons are either $B(\epsilon_N, \epsilon'_N)$ corners
or $B(\epsilon_S, \epsilon'_S)$ corners. This condition forces one
of the three white corners at $v$ to be either a white $(\epsilon_N,
\epsilon'_N)$ corner or a white $(\epsilon_S, \epsilon'_S)$ corner,
which is impossible by Lemma \ref{nowhitecorner}.

Therefore, every vertex in $\overline{G}$ has valency at least 4.
\end{proof}
Let $V, E, F$ be the number of vertices, edges and faces
respectively in $\overline{G}$. Since no trigons exist in $G_1$, no
trigons exist in $\overline{G}$. So every face except one in
$\overline{G}$ has at least 4 edges. Hence, 4$F- 4$ $\leq$ 2$E$.
Combined with $V - E + F = 2$, we get 2$E$ $\leq$ 4$V - 4$. On the
other hand, by the Claim, 4$V$ $\leq$ 2$E$. This is a contradiction.
\end{proof}

For the purpose of the next theorem, let $D^2(p, q)$ be a Seifert
fibered space over the disk with two exceptional fibers of orders
$p$ and $q$ respectively. Recall that $M(r_2)$ is a toroidal Dehn filling.
\begin{theorem} \label{23}
$M(r_2)$ is either $D^2(2, q)$ $\cup_T$ $D^2(2, s)$, $D^2(2, 3)$
$\cup_T$ $D^2(r, s)$ or $D^2(2, q)$ $\cup_T$ $D^2(3, s)$ where $T$
is an incompressible torus in $M(r_2)$. Furthermore, the fibers of
the two Seifert fibered subspaces of $M(r_2)$ intersect exactly once
on $T$.
\end{theorem}
\begin{proof}
$G_1$ contains a $W(\alpha, \delta)$
bigon and a $B(\alpha, \beta)$ face. Then Theorem \ref{Se} and Lemma \ref{Se1} imply that $M(r_2)$
is $D^2(2, q)$ $\cup_T$ $D^2(r, s)$ where $T$ is $\widehat{F}_2$ and
the two fibers are induced from an $\alpha \cup \delta$ curve and an
$\alpha \cup \beta$ curve respectively which intersect exactly once on
$T$.

On the other hand, Proposition \ref{trigon} guarantees the existence
of either a black bigon with two positive edges or a trigon in
$G_1$. If $G_1$ contains a black bigon with two positive edges,
$M(r_2)$ is $D^2(2, q)$ $\cup_T$ $D^2(2, s)$, which is the first
case.

If $G_1$ contains a trigon with three positive edges, which is a
good face, then by Theorem \ref{Se} one of orders of the four
exceptional fibers is 3. Then $M(r_2)$ is either $D^2(2, 3)$
$\cup_T$ $D^2(r, s)$ or $D^2(2, q)$ $\cup_T$ $D^2(3, s)$ depending
on which side the trigon lies.

The only remaining case is when $G_1$ has a trigon with two negative
edges and one positive edge. Let $f$ be a trigon with one
$\lambda$-edge and two negative edges where $\lambda \in \{\alpha,
\beta, \gamma, \delta\}$ as in Figure 28 ($i$). Then Lemma
\ref{order} (if $f$ is black) and Lemma \ref{nowhitecorner} (if $f$
is white) imply that the ($\epsilon, \epsilon'$) corner of $f$ must
be either a ($\epsilon_N, \epsilon'_S$) corner or a ($\epsilon_S,
\epsilon'_N$) corner.

Recall that $\widehat{F}_2$ is an incompressible torus containing
the graph $G_2$. Let $\widehat{H}$ be the intersection of $V_2$ and
the side of $M(r_2)$ where $f$ lies. $\partial f$ appears on $\widehat{F}_2 \cup
\partial \widehat{H}$ as in Figure 28 ($ii$). Let $A$ be an annulus
on $\widehat{F}_2$ containing all the edges of $f$. Let $N$ be the
regular neighborhood of $\widehat{H} \cup A \cup f$ in the side
where $f$ lies. Then it is easy to see that the fundamental group of
$N$ is isomorphic to the group $\langle$ $x, y \mid x^2y^3 = 1
\rangle$. Therefore $N$ is homeomorphic to $D^2(2, 3)$. This implies
the second case of $M(r_2)$.
\end{proof}
\begin{figure}[tbh]
\includegraphics{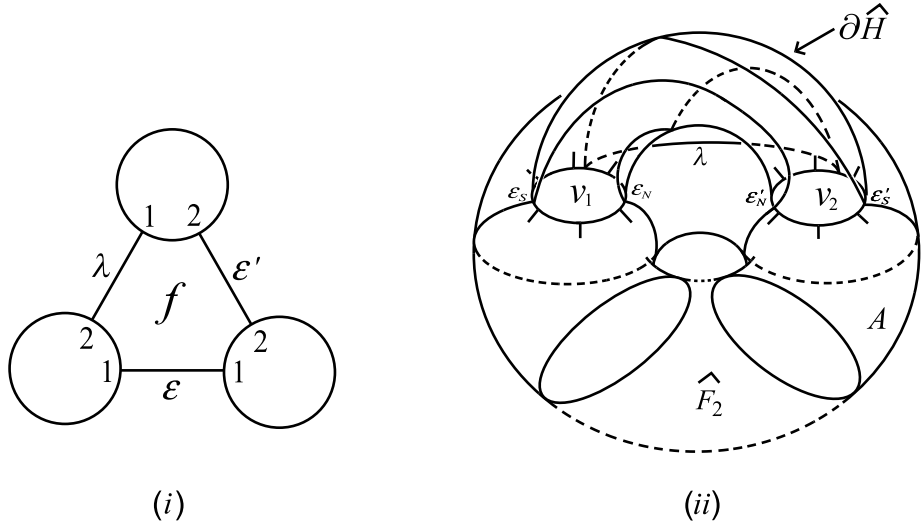}\caption{A trigon with two negative
edges and one positive edge.}\label{f:jj}
\end{figure}

\section{A link surgery description of the toroidal Dehn filling}
In this section, we interpret the toroidal Dehn filling $M(r_2)$ as
the manifold obtained by Dehn surgery on some link $L$ in $S^3$ by
using a black good face and a white good face.

$G_1$ contains a $B(\alpha, \beta)$ face and a $W(\alpha, \delta)$
bigon. We denote them by $f_1$ and $f_2$ respectively. Note that a $B(\alpha,
\beta)$ face is good by Lemma \ref{abface}.

Then we are exactly in the same situation as in \cite[Section
5]{GL4}. Thus we go through the same terminology and arguments as
in \cite[Section 5]{GL4}.

\begin{figure}[tbh]
\includegraphics{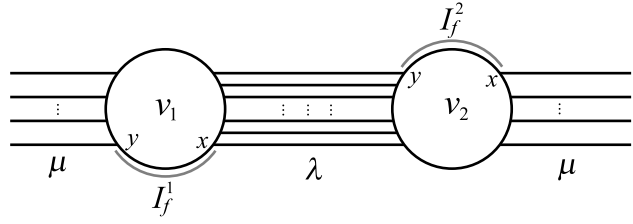}\caption{The external
interval of $f$.}\label{f:jj}
\end{figure}

\begin{figure}[tbh]
\includegraphics{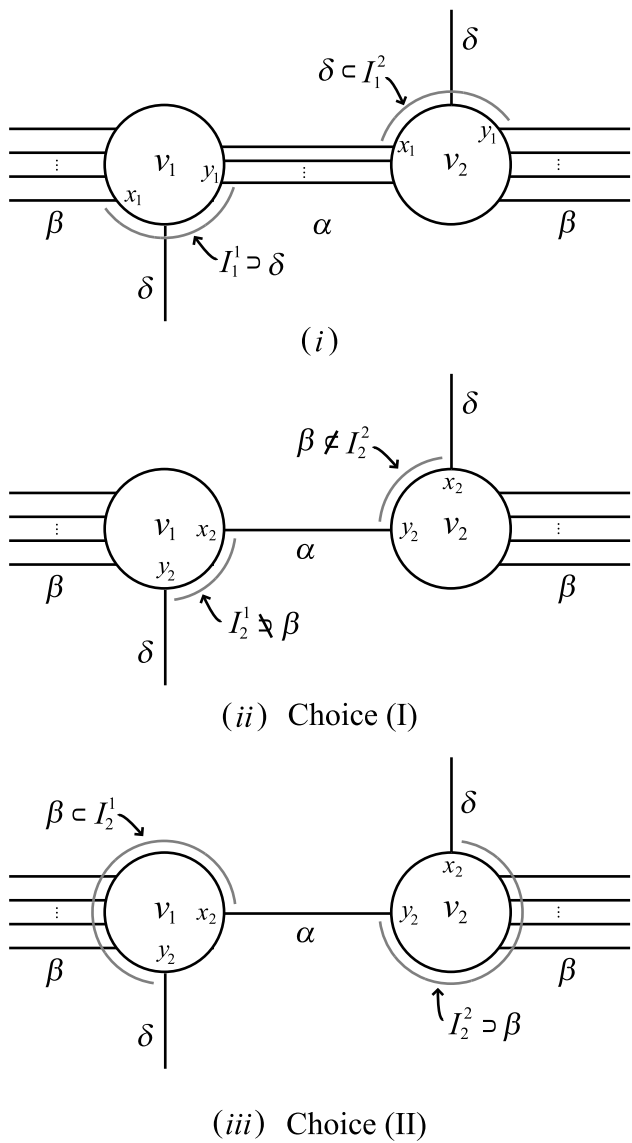}\caption{The external
intervals of $f_1$ and $f_2$.}\label{f:jj}
\end{figure}

Let $f$ be a good $(\lambda, \mu)$ face in $\Lambda$. Then by Lemma
4.2 in \cite{GL4}, there are vertices $x, y$ of $f$ in $\Lambda$,
representing corners of $\partial f$, such that the edges of
$\partial f$ on $G_2$ are, up to homeomorphism, as in Figure 29. The
interval [$x, y$] on vertex $v_i$, $i$ = 1, 2, that contains no edges
of $\partial f$ in its interior, is called \textit{the external
interval of f on vertex i}, denoted by $I^i_f$. In particular, if
$f$ is a bigon, there are two such intervals. Then we denote the
external interval of $f_i$ at vertex $v_j$ by $I^j_i$, $j= 1, 2$.
There are two choices on $I^j_2$:
\begin{itemize}
\item Choice (I). $I^j_2$ doesn't contain the endpoints of
the edge class $\beta$.
\item Choice (II). $I^j_2$ contains the endpoints of
the edge class $\beta$.
\end{itemize}

Note that $I^j_1$ contains the endpoints of the edge class $\delta$.
See Figure 30.

Recall that $M(r_2)$ = $\widehat{M}_B \cup_{\widehat{F}_2}
\widehat{M}_W$ and $\widehat{H}_B$ (resp. $\widehat{H}_W$) = $V_2$
$\cap$ $\widehat{M}_B$ (resp. $\widehat{M}_W$). For convenience, we
change indices $B, W$ into $1, 2$ respectively and we let $V$ be the
attached solid torus to $M$. Recall that $r_1$ and $r_2$ are slopes on
$\partial_0M=\partial V$ such that $M(r_1)$ is reducible and
$M(r_2)$ is toroidal and the boundary components of $F_1$ have
slopes $r_1$.

Now follow the argument from the page 452 through the page 457 in
\cite[Section 5]{GL4} with some notations changed, i.e. with
$\epsilon,$ $\delta_1,$ $\delta_2,$ $H_1,$ $H_2,$ $\mu,$ $\tau,$
$Q,$ $T$ and $M_i$ replaced by $\alpha,$ $\beta,$ $\delta,$
$\widehat{H}_1,$ $\widehat{H}_2,$ $r_1,$ $r_2,$ $F_1,$ $F_2,$ and
$\widehat{M}_i$.

\begin{figure}[tbh]
\includegraphics{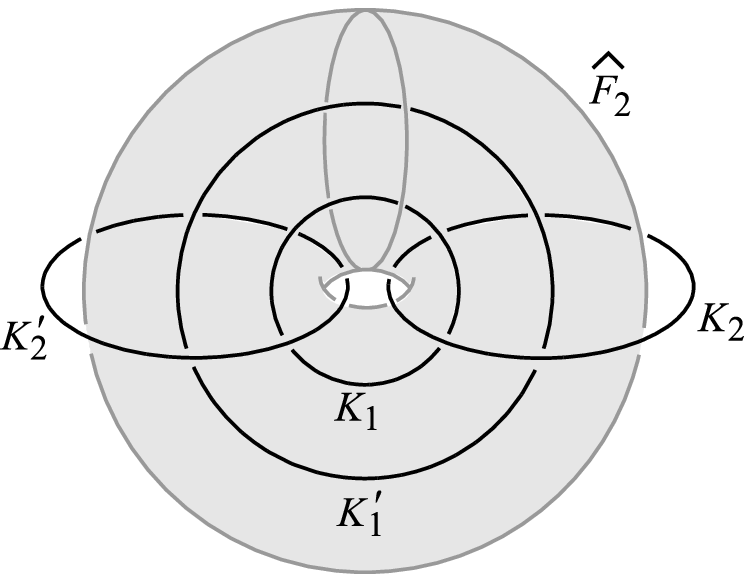}\caption{The link $L_0=K_1 \cup K'_1 \cup K_2 \cup K'_2$.}\label{f:jj}
\end{figure}

\begin{figure}[tbh]
\includegraphics{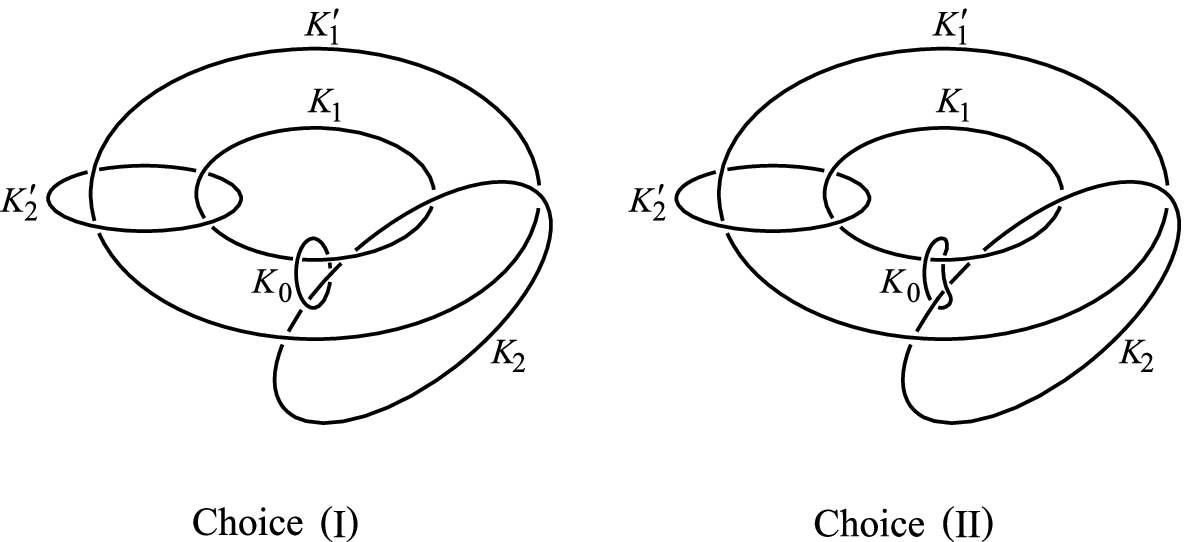}\caption{The link $L=K_1 \cup K'_1 \cup K_2 \cup K'_2 \cup K_0$.}\label{f:jj}
\end{figure}

\begin{figure}[tbh]
\includegraphics{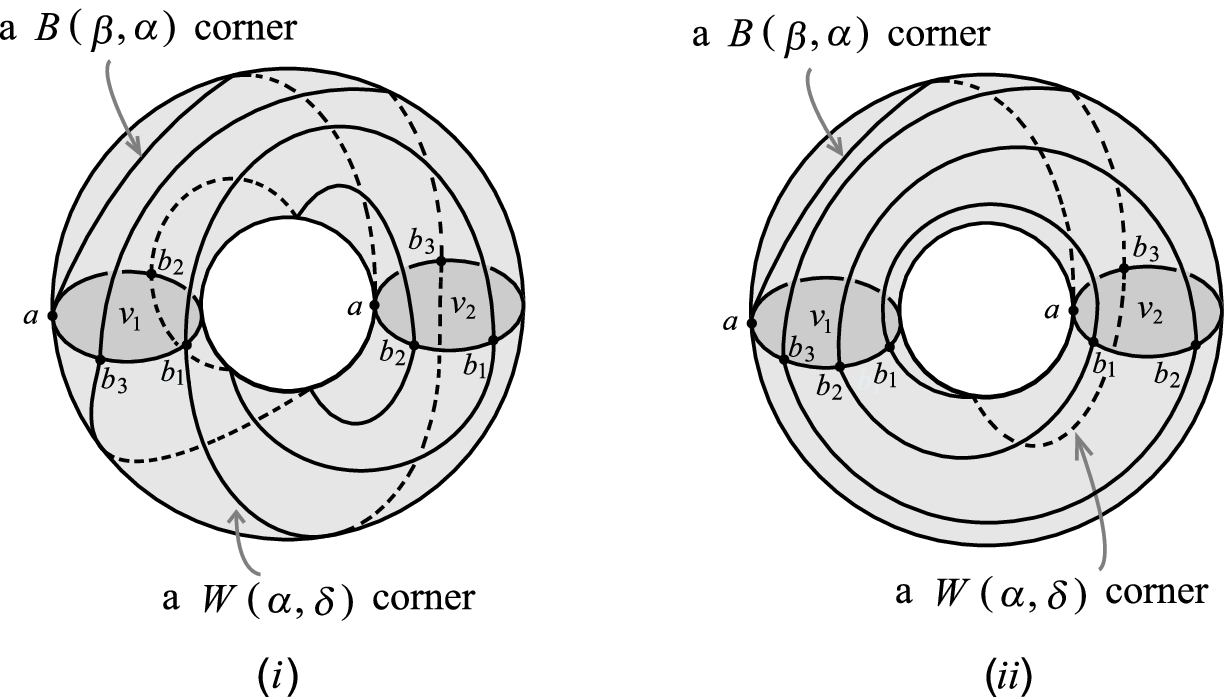}\caption{The $r_1$-framing.}\label{f:jj}
\end{figure}

Then from the page 456 and the page 457 in \cite[Section 5]{GL4}, we have the
following: \textit{$(M(r_2), K_{r_2})$ is obtained from $(S^3, K_0)$
by some Dehn surgery on $L_0$ where $K_{r_2}$ is the core of $V$ in
$M(r_2)$. That is, the exterior of $K_{r_2}$ is obtained from the
exterior of $K_0$ in $S^3$ by a Dehn surgery on $L_0$. Moreover, given
the $0$-framing of $K_0$ in $S^3$, the meridian of $K_0$ corresponds to the
slope $r_2$ on $\partial V$.}

Here, as in Section 5 in \cite{GL4}, $K_0$ means the core of $V =
\widehat{H}_1 \cup \widehat{H}_2 \subset S^3$ and $L_0$ means the
link $K_1 \cup K'_1 \cup K_2 \cup K'_2$ in $S^3$ depicted in Figure 31. Also
let $L$ be the link $L = K_1 \cup K'_1 \cup K_2 \cup K'_2 \cup K_0$.
Then $L$ looks like Figure 32 according to the two choices (I), (II)
on $I^j_2$. Note that $L$ in Choice (II) is the reflection of $L$ in
Choice (I) (in which $K_i \mapsto K_i$, $0 \leq i \leq 2$, and $K'_i
\mapsto K'_i$, $i= 1, 2$.). Therefore under this reflection the
$r_1$-framing on $K_0$ goes to the ($-r_1$)-framing (and the
longitude to the longitude).

Let $L(\theta, \varphi, \psi, \omega, \pi)$ be the Dehn filling on
the exterior of $L$ where $\theta, \varphi, \psi, \omega, \pi$
denote the filling slopes on $K_1, K'_1, K_2, K'_2, K_0$
respectively. When one of the filling slopes is an asterisk
($\ast$), then no filling is done on the corresponding component.

\begin{proposition} \label{slope1}
$M = L(\theta, \varphi, \psi, \omega, \ast)$ where $L$ is the link
from Choice $\rm(I)$ as in Figure $32$. Furthermore, the reducible
Dehn filling $M(r_1)$ is $L(\theta, \varphi,$ $\psi, \omega, 1/3)$
and the toroidal Dehn filling $M(r_2)$ is $L(\theta,$ $\varphi,$
$\psi,$ $\omega,$ $1/0)$.
\end{proposition}
\begin{proof}
>From the above discussion, we need only to show that $\pi$ = 1/3 for $M(r_1)$.
In other words, with respect to the
framing of $K_0$, the slope $r_1$ is 1/3. Let $a$ be a vertex in
$G_1$ which contains a $B(\beta, \alpha)$ corner of a $B(\alpha,
\beta)$ face $f_1$. That is, $v_1$ in $G_2$ has the label $a$ at the endpoint
of a $\beta$-edge and $v_2$ in $G_2$ has the label $a$ at the endpoint of an
$\alpha$-edge. Let $b$ be a vertex in $G_1$ which contains a $W(\alpha,
\delta)$ corner of a $W(\alpha, \delta)$ bigon $f_2$. Then it
suffices to show that $\partial b$ has slope 1/3.

$\partial b$ consists of three black corners and three white
corners, one of which is a $W(\alpha, \delta)$ corner. There are
three occurrences of the label $b$ at $v_1$ and $v_2$. Let $b_1, b_2$
and $b_3$ be the labels representing three occurrences of $b$ at $v_1$
and $v_2$ such that $b_i$ at $v_1$ and $b_i$ at $v_2$ are connected
by a black corner of the vertex $b$, and $b_i$ at $v_2$ and
$b_{i+1}$ at $v_1$ are connected by a white corner of the vertex
$b$, $i$ = 1, 2, 3 (mod 4). We assume that a $W(\alpha, \delta)$
corner connects $b_3$ at $v_2$ and $b_1$ at $v_1$. The labels $b_1,
b_2$ and $b_3$ on $\partial v_1$ and $\partial v_2$ occur either in the order $b_1
b_2b_3$ or in the order $b_1b_3b_2$.

First, we assume the order $b_1 b_2b_3$. Then we take Choice (I) for
$I^j_2$. Since the $r_1$-framings on $\widehat{H}_1$ and
$\widehat{H}_2$ are shown Figure 5.16 and Figure 5.17
Case (I) in \cite{GL4} respectively, a $B(\beta, \alpha)$ corner and a $W(\alpha,
\delta)$ corner, which are uniquely determined on $\partial
\widehat{H}_1$ and $\partial \widehat{H}_2$ respectively are given
in Figure 33 $(i)$. Therefore, $\partial v$ has slope 1/3.

Second, we assume the order $b_1b_3b_2$. Then we take Choice (II)
for $I^j_2$. The $r_1$-framing on $\widehat{H}_2$ is shown in Figure
5.17 Case (II) in \cite{GL4}. Thus, $\partial v$ has slope $-1/3$. See Figure 33
($ii$). Under the reflection sending $L$ in Choice (II) to $L$ in
Choice (I), $-1/3$ in Choice (II) goes to 1/3 in Choice (I).
Therefore we may assume without loss of generality that $\pi=1/3$ in
Choice (I). This completes the proof of $\pi= 1/3$.
\end{proof}

\begin{lemma} \label{psi}
$\psi$ is equal to $2$ in Choice $\rm(I)$ and $-2$ in Choice $\rm(II)$.
\end{lemma}
\begin{proof}
$\psi$ is the filling slope on $K_2$, which is a core of $V_2$
disjoint from $\widehat{H}_2 \cup R_2$. Also ($W_2; \widehat{H}_2,
R_2$) is obtained from ($V_2; \widehat{H}_2, R_2$) by $\psi$-Dehn
filling on $K_2$ as described in the page 456 in \cite{GL4}. Recall that $W_2$ = nbhd($A_2 \cup \widehat{H}_2
\cup f_2$). Since $f_2$ is a $W(\alpha, \delta)$ bigon, $W_2$ is the
neighborhood of a M$\ddot{\rm{o}}$bius band. Then with the standard framing on
$K_2$ which is shown in Figure 5.17 Case (I) in \cite{GL4}, we easily see that $\psi$ =
2 in Choice (I). For Choice (II), notice that the longitude has the
opposite orientation to Choice (I). Therefore, $\psi$ = $-2$ in
Choice (II).
\end{proof}

\section{Dehn fillings on a link and Tangle fillings on a tangle}

>From now on, $L$ is the link in Choice (I) as in Figure 32.
Proposition \ref{slope1} says that $M$, the reducible Dehn filling
$M(r_1)$ and the toroidal Dehn filling $M(r_2)$ are obtained by Dehn
surgery on $L$. $L$ is strongly invertible and the quotient of
the exterior under this involution gives rise to a
tangle. Therefore $M$, $M(r_1)$ and $M(r_2)$ are double branched
cover of $S^3$ along some tanlges or some links which are obtained
by rational tangle fillings on the tangle obtained from the quotient of
the exterior of $L$.
\begin{figure}[tbh]
\includegraphics{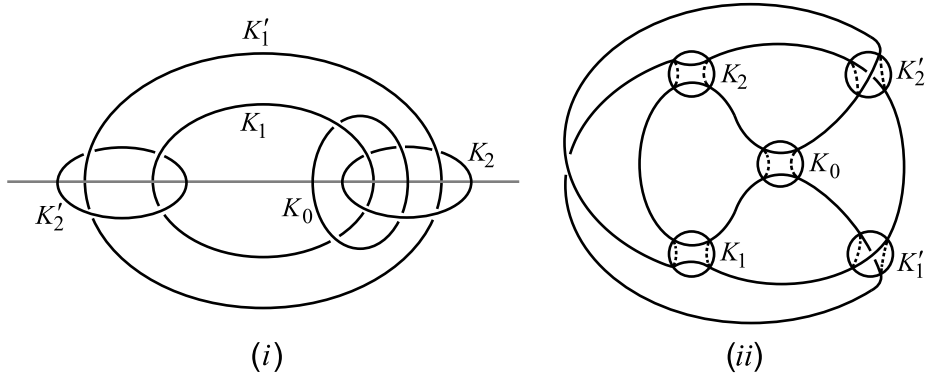}\caption{The strong invertibility of $L$ and the quotient of ($S^3, L$).}\label{f:jj}
\end{figure}

We will use the same convention on rational tangles and the same
terminology in \cite[Section 6]{GL4}. Let $T(\alpha_1,
\cdot\cdot\cdot, \alpha_n)$ be a filling of a tangle $T$ where
$\alpha_i \in  \bf{Q}$ $\cup$ $\{1/0\}$, $1 \leq i \leq n$. Then the
double branched cover of $S^3$ along the link $T(\alpha_1,
\cdot\cdot\cdot, \alpha_n)$ is a closed 3-manifold denoted by
$\widetilde{T}(\alpha_1, \cdot\cdot\cdot, \alpha_n)$.

By Proposition \ref{slope1}, $M$ = $L(\theta, \varphi, \psi, \omega,
\ast)$. Figure 34 ($i$) shows that $L$ is strongly invertible and
depicts the quotient along with its fixed set. Then the quotient of
($S^3, L$) under this strong inversion is $P(\frac{0}{1},
\frac{1}{1}, \frac{0}{1}, \frac{1}{1}, \frac{0}{1})$ as in Figure 34
$(ii)$, where $P$ is the pentangle (See \cite[Section 6]{GL4} for the
definition of the pentangle.) and invariant regular neighborhoods of
$L$ become tangle balls in the quotient. In other words, the
quotient of the exterior $E(L)$ of $L$ in $S^3$ under
the involution is the pentangle $P$. The longitude framings of
$\partial E(L)$ pass to the $\frac{1}{0}$ framings of the tangle
balls as indicated by the dashed line. The meridian framings of
$\partial E(L)$ pass to the $\frac{0}{1}$ framings of the tangle
balls for $K_0, K_1, K_2$ and $\frac{1}{1}$ framings for $K'_1,
K'_2$ as indicated by the solid line. In particular the (meridian,
longitude)-framings of $K_0, K_2$ pass to the ($\frac{0}{1},
\frac{1}{0}$)-framing of the corresponding tangle spheres of $P$,
which implies with
Proposition \ref{slope1} and Lemma \ref{psi} that the
reducible Dehn filling $M(r_1)$ = $L(\theta, \varphi, 2, \omega,
1/3)$ is $\widetilde{P}(\theta', \varphi', 1/2, \omega', 3)$ and the
toroidal Dehn filling $M(r_2)$ = $L(\theta, \varphi, 2, \omega,
1/0)$ is $\widetilde{P}(\theta', \varphi', 1/2, \omega', 0/1)$.

\begin{figure}[tbh]
\includegraphics{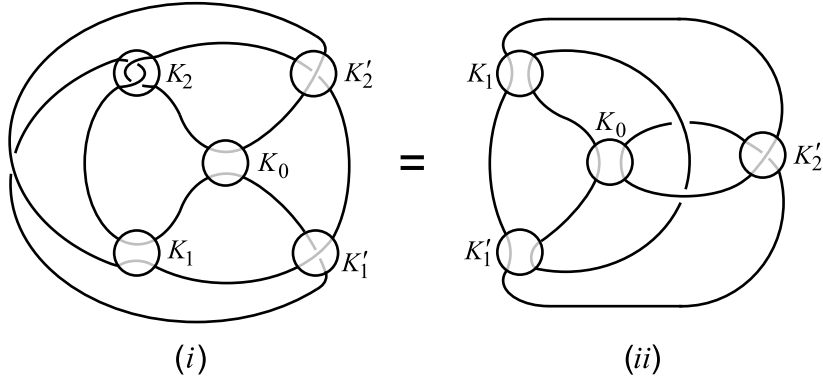}\caption{The pentangle $P(\theta', \varphi', 1/2,
\omega', \ast)$ = the tangle $Q(-1/\theta', -1/\varphi',
-1/\omega', \ast)$.}\label{f:jj}
\end{figure}

However we can regard the pentangle $P(\theta', \varphi', 1/2,
\omega', \ast)$ as the tangle $Q(-1/\theta', -1/\varphi',
-1/\omega', \ast)$ by rotating the pentangle by $-\pi$/2 where $Q$
is the tangle as shown in Figure 35 ($ii$). Note that the tangle $Q$
admits a symmetry by rotation in the horizontal axis. Without
ambiguity, if we let $\theta, \varphi, \omega$ be $-1/\theta'$,
$-1/\varphi'$, $-1/\omega'$ respectively, then since $M$ is
$\widetilde{P}(\theta', \varphi', 1/2, \omega', \ast)$, $M$ =
$\widetilde{Q}( \theta, \varphi, \omega, \ast)$, $M(r_1)$ =
$\widetilde{Q}( \theta, \varphi, \omega, -1/3)$ and $M(r_2)$ =
$\widetilde{Q}$ $( \theta, \varphi, \omega, 1/0)$.

\begin{theorem} \label{main2}
$M$ is the double branched cover of one of the tangles $Q( 2, p-2,
1/p, \ast)$ where $p$ is an integer and $p \leq -3$.
\end{theorem}

To prove this, we do a series of computations describing the double
branched covers of certain tangle filling of $Q$ and combining this
with certain Dehn filling theorems, we put restrictions on the
tangle slopes $\theta, \varphi$ and $\omega$.
\begin{lemma} \label{rr}
Let $M$ be an irreducible $3$-manifold with $\partial M$ a torus.
Let $\theta, \varphi$ be slopes on $\partial M$ such that
$M(\theta)$ is reducible and $M(\varphi)$ is either reducible, a
lens space, or $S^3$. Then $\Delta( \theta, \varphi) \leq 1$.
\end{lemma}
\begin{proof}
This is exactly Lemma 8.1 of \cite{GL4}.
\end{proof}

\begin{figure}[tbh]
\includegraphics{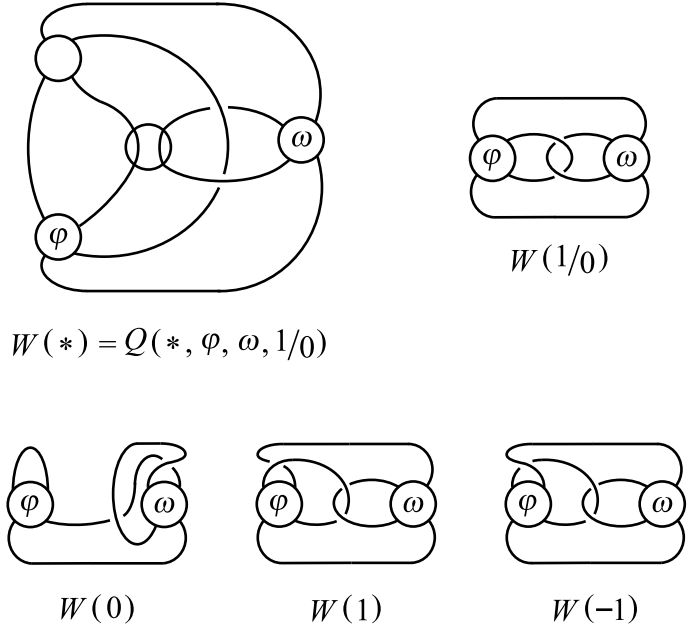}\caption{$W(\ast)=Q( \ast, \varphi, \omega, 1/0)$ and the rational tangle fillings.}\label{f:jj}
\end{figure}

\begin{figure}[tbh]
\includegraphics{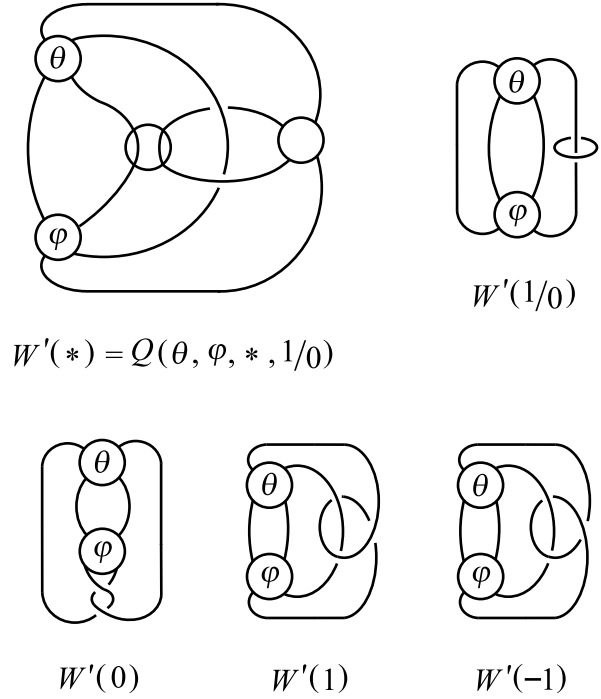}\caption{$W'(\ast)=Q( \theta, \varphi, \ast, 1/0)$ and the rational tangle fillings.}\label{f:jj}
\end{figure}

\begin{lemma} \label{zeroone}
$\theta, \varphi \not\in \{ -1, 0, 1, 1/0\}$ and $\omega \not\in \{
-1, 0, 1, 1/0, 1/2\}$.
\end{lemma}
\begin{proof}
Let $W(\ast)$ be $Q( \ast, \varphi, \omega, 1/0)$. Then $W(\theta)$
= $Q(\theta, \varphi, \omega, 1/0)$ and $\widetilde{W}(\theta)$ =
$\widetilde{Q}(\theta, \varphi, \omega, 1/0)$( = $M(r_2)$). Theorem
\ref{23} says that $M(r_2)=\widetilde{W}(\theta)$ is a toroidal manifold
and is not a Seifert fibered space. Moreover, by Lemma \ref{rr},
$\widetilde{W}(\theta)$ is also irreducible since $\Delta(r_1, r_2)=
\Delta( 1/0, -1/3)= 3$. Consider the rational tangle fillings
$W(-1), W(0), W(1)$ and $W(1/0)$. See Figure 36. Then the double
branched cover of each tangle filling is either a Seifert fibered
space or a reducible manifold, which is impossible. Therefore,
$\theta \not\in \{ -1, 0, 1, 1/0\}$. On the other hand, the
tangle $Q$ admits a symmetry by rotation in the horizontal axis,
which interchanges $\theta, \varphi$. Therefore $ \varphi \not\in \{ -1,
0, 1, 1/0\}$.

For $\omega$, let $W'(\ast)$ be $Q( \theta, \varphi, \ast, 1/0)$.
Then $M(r_2)$ = $\widetilde{W'}(\omega)$ is a toroidal and
irreducible manifold and is not a Seifert fibered space. Consider
the rational tangle fillings illustrated in Figure 37. The double
branched cover of each tangle filling is either a Seifert fibered
space or a reducible manifold, which is impossible. Hence $ \omega \not\in \{ -1,
0, 1, 1/0\}$.

If $\omega=1/2$, then $M(r_2)=\widetilde{W}'(1/2)$ contains a Klein bottle from
figure 37, which is a contradiction. Thus
$\omega \neq 1/2$.
\end{proof}

\begin{figure}[tbh]
\includegraphics{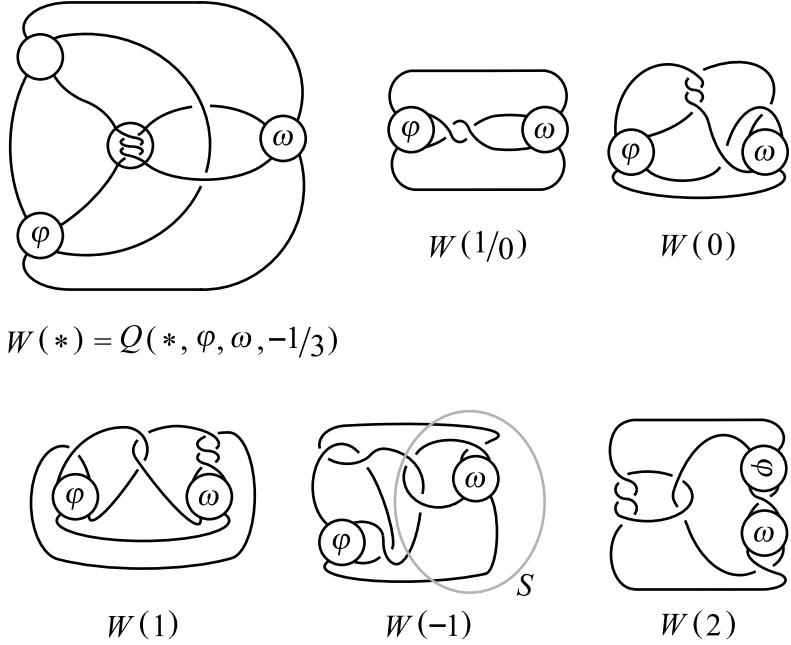}\caption{$W(\ast)=Q( \ast, \varphi, \omega, -1/3)$ and the rational tangle fillings.}\label{f:jj}
\end{figure}

\begin{figure}[tbh]
\includegraphics{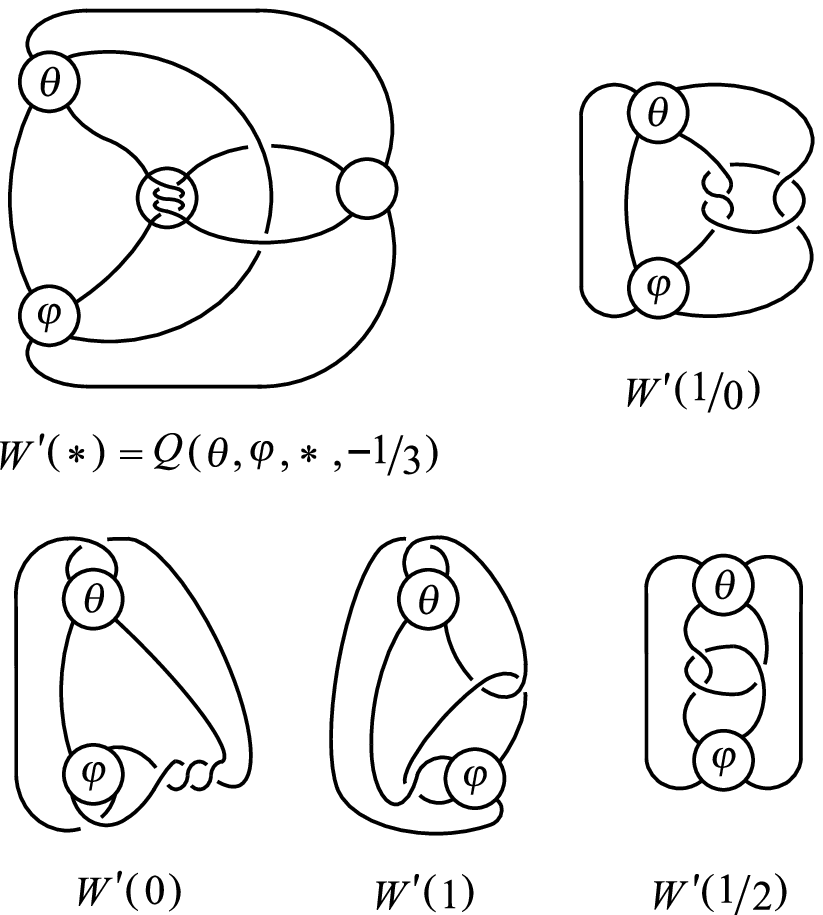}\caption{$W'(\ast) = Q( \theta, \varphi,
\ast, -1/3)$ and the rational tangle fillings.}\label{f:jj}
\end{figure}

Recall that $D^2(a, b)$ is a Seifert fibered space over $D^2$ with
two exceptional fibers of orders $a, b$. Similarly, let $S^2(a, b,
c)$ be a Seifert fibered space over $S^2$ with three exceptional
fibers of orders $a, b, c$.
\begin{lemma} \label{integers}
$\theta, \varphi$ are integers and $\omega$ is 1/$k$ where $k$ is an
integer.
\end{lemma}
\begin{proof}
First, we show that $\theta$ is an integer. Let $W(\ast) = Q( \ast,
\varphi, \omega, -1/3)$. Then $\widetilde{W}(\theta) =
\widetilde{Q}( \theta, \varphi, \omega, -1/3)$ is $M(r_1)$, which is
reducible. We perform the rational tangle fillings shown in Figure
38. Observe that $\widetilde{W}(1/0)$ is either $S^3$, $S^2 \times
S^1$ or a lens space. Hence $\Delta(\theta, 1/0) \leq 1$ by Lemma
\ref{rr}, provided that $\widetilde{W}(\ast)$ is irreducible. Since
$\theta \neq 1/0$ by Lemma \ref{zeroone}, $\theta$ is an integer.
Thus, to show that $\theta$ is an integer, it remains to show that
$\widetilde{W}(\ast)$ is irreducible.

Suppose for contradiction that $\widetilde{W}(\ast)$ is reducible,
i.e. $\widetilde{W}(\ast) = X$ $\#$ $Y(\ast)$ where $Y$ contains
the torus boundary component where Dehn fillings are performed. Note
that if there are two slopes $p$ and $q$ such that
$\widetilde{W}(p)$ and $\widetilde{W}(q)$ are prime manifolds, then
$\widetilde{W}(p)$ = $\widetilde{W}(q)$. However
$\widetilde{W}(1/0)$ is either $S^3$, $S^2 \times S^1$ or a lens
space, which are all prime. Thus $\widetilde{W}(p)$ is either $S^3$,
$S^2 \times S^1$ or a lens space if $\widetilde{W}(p)$ is prime.
>From Figure 38, $\widetilde{W}(0)$ and $\widetilde{W}(-1)$ are prime
manifolds because $\varphi, \omega \not\in \{ -1, 0, 1, 1/0\}$,
whence they must be either $S^3$, $S^2 \times S^1$ or a lens space.

Suppose $\varphi \neq 1+1/a$ for an integer $a$. If $\omega \neq
1/(1+1/b)$ for an integer $b$, then $\widetilde{W}(-1)$ is a
toroidal prime manifold, which is neither $S^3, S^2 \times S^1$ nor
a lens space. If $\omega = 1/(1+1/b)$ for some integer $b$, then
$\widetilde{W}(-1)$ is $S^2(e, f, g)$ where $e, f, g \geq 2$ since
$\omega \neq 0, 1/2$, which is neither $S^3, S^2 \times S^1$ nor a
lens space, a contradiction. Therefore $\varphi$ must be 1+1/$a$ for
some integer $a$. Then $\widetilde{W}(1)$ = $S^2(|2a+1|, 2, t)$ for
some integer $t$. Note that $|2a+1|$ $\neq 0, 1$ because $\varphi
\not\in \{ -1, 0, 1, 1/0\}$. If $\omega$ $\neq 1/3$ i.e. $t$ $\neq
0$, then $\widetilde{W}(1)$ is a prime manifold, which must be either
$S^3, S^2 \times S^1$ or a lens space. Hence $\omega$ must be
1/(3+1/$s$) for some integer $s$. This implies that
$\widetilde{W}(-1)$ = $S^2(|2a+1|, 2, |2s+1|)$. However, since
$\varphi$ $\neq -1, 0, 1$ and $\omega$ $\neq -1, 0, 1, 1/2$,
$|2a+1|$ $\geq 2$ and $|2s+1|$ $\geq 2$, whence $\widetilde{W}(-1)$
is neither $S^3, S^2 \times S^1$ nor a lens space, a contradiction.
If $\omega$ = $1/3$, then $\widetilde{W}(-1)$ = $S^2(|2a+1|, 2, 2)$,
which is neither $S^3, S^2 \times S^1$ nor a lens space, a
contradiction. This completes the proof that $\widetilde{W}(\ast)$
is irreducible.

By the symmetry by rotation in the horizontal axis on the tangle
$Q$, $\varphi$ is also an integer.

To show that $\omega$ is 1/$k$ for some integer $k$, a similar
argument to the above applies. Let $W'(\ast) = Q( \theta, \varphi,
\ast, -1/3)$. Then $\widetilde{W'}(\omega) = \widetilde{Q}( \theta,
\varphi, \omega, -1/3)$ is $M(r_1)$, which is reducible. We perform
the rational tangle fillings shown in Figure 39. Since
$\widetilde{W'}(0)$ is either $S^3, S^2 \times S^1$ or a lens space,
$\omega$ is 1/$k$ for some integer $k$ by Lemma \ref{rr}, provided
that $\widetilde{W'}(\ast)$ is irreducible. However,
$\widetilde{W'}(1/2)$ is prime since $\theta, \varphi \not\in \{ -1,
0, 1, 1/0\}$, and $\widetilde{W'}(1/2)$ $\neq$ $\widetilde{W'}(0)$.
Thus, $\widetilde{W'}(\ast)$ is irreducible. The proof is complete.
\end{proof}
\begin{proposition} \label{phi}
$\theta$ is $2$ up to the symmetry by rotation in the horizontal
axis on $Q$.
\end{proposition}

Before proving Proposition \ref{phi}, we need to consider several
special cases.
\begin{lemma} \label{-2}
$\theta$, $\varphi \neq$ $-2$.
\end{lemma}
\begin{proof}
Suppose $\varphi=-2$. Let $W(\ast) = Q(\ast, -2, \omega, -1/3)$. Then
$\widetilde{W}(\theta)$ is reducible. Consider the rational tangle fillings
shown in Figure 38 with the $-2$-rational tangle substituted for
$\varphi$. Note that $\widetilde{W}(1)$ is a lens space or a
reducible manifold. Since $\widetilde{W}(1/0)$ and
$\widetilde{W}(-1)$ are prime and they are not homeomorphic,
$\widetilde{W}(\ast)$ is irreducible. Therefore we can apply Lemma
\ref{rr} to get $\Delta(\theta, 1)$ $\leq$ 1. Thus $\theta$ = 0, 1,
2. However 0 and 1 are impossible by Lemma \ref{zeroone}. Then
$\theta$ = 2. We will show that 2 is also impossible.

Consider the toroidal Dehn filling $M(r_2)$ = $\widetilde{Q}(\theta,
\varphi, \omega, 1/0)$. Since $\theta=2$ and $\varphi=-2$, $M(r_2)$
contains a Klein bottle, which is a contradiction to the assumption
that $M(r_2)$ does not contain a Klein bottle.

By the symmetry on $Q$, $\theta \neq -2$.
\end{proof}

\begin{lemma} \label{-3}
If $\varphi =$ $-3$, then $\theta =$ $2$.
\end{lemma}
\begin{proof}
Assume that $\varphi= -3$. Let $W(\ast) = Q(\ast, -3, \omega,
-1/3)$. Then $\widetilde{W}(\theta)$ is reducible. Consider the
rational tangle fillings shown in Figure 38 with the $-3$-rational tangle
substituted for $\varphi$.

Assume $\widetilde{W}(\ast)$ is hyperbolic. Notice that
$\widetilde{W}(1)$ = $S^2(2, 2, |1/\omega-3|)$ and
$\widetilde{W}(-1)$ is toroidal since $\omega \not\in \{ -1, 0, 1,
1/0, 1/2\}$. Then $\widetilde{W}(\ast)$ admits a reducible Dehn
filling $\widetilde{W}(\theta)$, a toroidal Dehn filling $\widetilde{W}(-1)$, and a Dehn filling of the form
$S^2(2, 2, 1/\omega-3)$. By Corollary 1.3 in \cite{Lee},
$\Delta(\theta, 1)$ $\leq$ 2 and by Theorem 1.1 in \cite{O1},
$\Delta(\theta, -1)$ $\leq$ 3. Thus $\theta$ = $-1$, 0, 1 or 2.
However, Lemma \ref{zeroone} leads to $\theta$ = 2. The proof is
complete under the assumption that $\widetilde{W}(\ast)$ is
hyperbolic.

We show that $\widetilde{W}(\ast)$ is hyperbolic. Equivalently, it
suffices to show that $\widetilde{W}(\ast)$ is irreducible,
$\partial$-irreducible, non-Seifert fibered and atoroidal. Since
$\widetilde{W}(1/0)$ and $\widetilde{W}(-1)$ are prime and are not
homeomorphic, $\widetilde{W}(\ast)$ is irreducible.

Suppose $\widetilde{W}(\ast)$ is a Seifert fibered space. Since
$\widetilde{W}(1/0)$ is a lens space, $\widetilde{W}(\ast)$ is a
Seifert fibered space over $D^2$ with at most 2 exceptional fibers.
Hence $\widetilde{W}(r)$ is a Seifert fibered space over $S^2$ with
at most three exceptional fibers except for one slope $r$, for which
$\widetilde{W}(r)$ is reducible. Figure 38 shows that
$\widetilde{W}(-1)$ is irreducible. Hence $\widetilde{W}(-1)$ = $S^2(e, f, g)$
where $e, f, g \geq 1$, which doesn't contain a separating
incompressible torus. On the other hand, $\widetilde{W}(-1)$ contains a separating
incompressible torus which is the double branched cover of $S$ in Figure 38, since
$\omega \not\in \{-1, 0, 1, 1/0, 1/2\}$. This is a contradiction. Therefore
$\widetilde{W}(\ast)$ is not a Seifert fibered space.

Suppose $\widetilde{W}(\ast)$ is $\partial$-irreducible. After
$\partial$-compression, the torus boundary becomes a sphere which
must bound a 3-ball since $\widetilde{W}(\ast)$ is irreducible. This
implies that $\widetilde{W}(\ast)$ is a solid torus, which is a
Seifert fibered space, a contradiction.

Suppose $\widetilde{W}(\ast)$ is toroidal. Let $T$ be an essential
torus in $\widetilde{W}(\ast)$. Since $\widetilde{W}(1/0)$ is a lens
space, $T$ must be separating, which allows us to let
$\widetilde{W}(\ast)$ = $A \cup_T B(\ast)$ where $B$ contains the
torus boundary of $\widetilde{W}(\ast)$.

\begin{claim}
$T$ is compressible in $\widetilde{W}(-1)$.
\end{claim}
\begin{proof}
Suppose $T$ is incompressible in $\widetilde{W}(-1)$. Figure 38
shows that there is a unique incompressible torus up to isotopy in
$\widetilde{W}(-1)$, which must be the double branched cover of the
sphere $S$ as in Figure 38. Therefore $\widetilde{W}(-1)$ = $D^2(2,
4)$ $\cup_T$ $D^2(2, |1/\omega-1|)$. Thus $A$ is either $D^2(2, 4)$
or $D^2(2, |1/\omega-1|)$. Observe from Figure 38 that
$\widetilde{W}(0)$ = $S^2(3, 3, |1/\omega-2|)$, $\widetilde{W}(1)$ =
$S^2(2, 2, |1/\omega-3|)$ and $\widetilde{W}(2)$ = $S^2(3, 2,
|1/\omega+1|)$. Hence $A$ is $D^2(2,3)$ and $\omega$ = 1/4. $A$ =
$D^2(2, 3)$ is the complement of the trefoil knot i.e. the $(2,3)$
torus knot $T_{2,3}$. $\widetilde{W}(1)$ = $S^2(2, 2, 1)$ becomes a
lens space. It is easy to see that $\widetilde{W}(1)$ = $L(8, -3)$.
Since $B(1)$ is a solid torus, $\widetilde{W}(1)$ = $L(8, -3)$ is
obtained by Dehn filling on $A$ = $T_{2,3}$. This is impossible by
Corollary 7.4 in \cite{G2}.
\end{proof}

$T$ is compressible in $\widetilde{W}(1/0)$, $\widetilde{W}(0)$,
$\widetilde{W}(1)$, $\widetilde{W}(-1)$ and $\widetilde{W}(2)$.
Since $\Delta(2, -1)$ = $\Delta(2, 0)$ = $\Delta(1, -1)$ $\geq 2$,
by Lemma 2.4 in \cite{EW} $B(\ast)$ is a cable space $C(p, q)$
with the cabling slope 1/0. (See Section 3 in \cite{GL} for the
definition of a cable space.). Moreover, since $\Delta(\theta, 1/0)$
= 1, $B(\theta)$ is a solid torus, which makes $T$ compressible in
$\widetilde{W}(\theta)$. Therefore $\widetilde{W}(1/0)$ =
$A(r_{1/0})$ $\#$ $L(q, p)$ and $\widetilde{W}(\theta)$ = $A(r_\theta)$ for some
slopes $r_{1/0}, r_\theta$ on $T$ with $\Delta(r_{1/0}, r_\theta)$ =
$|q|\Delta(\theta, 1/0)$ = $|q|$. Since $\widetilde{W}(1/0)$ is a lens space,
$A(r_{1/0})$ must be $S^3$. Then $|q|$ must be 1 by Lemma
\ref{rr} since $A$ is irreducible, $A(r_{1/0})$ is $S^3$ and
$A(r_\theta)$ is reducible. The fact that $q$ = 1 or $-1$ implies
that $T$ is $\partial$-parallel in $\widetilde{W}(\ast)$, which is a
contradiction to the assumption that $T$ is essential in
$\widetilde{W}(\ast)$. Thus $\widetilde{W}(\ast)$ is atoroidal.
\end{proof}

\begin{lemma} \label{3}
If $\varphi =$ $3$, then $\theta =$ $2$.
\end{lemma}
\begin{proof}
Assume that $\varphi=3$. Let $W(\ast) = Q(\ast, 3, \omega, -1/3)$.
Then $\widetilde{W}(\theta)$ is reducible. Consider the rational
tangle fillings shown in Figure 38 with the 3-rational tangle
substituted for $\varphi$.

Assume $\widetilde{W}(\ast)$ is hyperbolic. $\widetilde{W}(-1)$ is
toroidal since $\omega \not\in \{ -1, 0, 1,$ $1/0, 1/2\}$.
$\widetilde{W}(\theta)$ is reducible. Thus $\Delta(\theta, -1)$
$\leq$ 3. The possible values of $\theta$ are $-4, -3, -2, 2$ since
$\theta \not\in \{ -1, 0, 1, 1/0\}$. If $\theta$ is 2, we are done.
If $\theta$ is $-3$, then this belongs to the case of $\varphi$ =
$-3$ by the symmetry on $Q$, whence we are done. Also, $\theta \neq -2$
by Lemma \ref{-2}. The remaining case is $\theta=-4$.

Suppose $\theta = -4$. Let $W'(\ast) = Q(-4, 3, \ast, -1/3)$. Then
$\widetilde{W}'(\omega)$ is reducible. Similarly, we can easily see
from Figure 39 that $\widetilde{W}'(1/0)$ is the lens space $L(7,
-2)$ and $\widetilde{W}'(1)$ is $S^2(4, 2, 3)$. They are
non-homeomorphic prime manifolds. Thus $\widetilde{W}'(\ast)$ is irreducible.
Moreover, Lemma \ref{rr} leads to $\Delta(1/0, \omega)$ $\leq$ 1. Since
$\omega$ is of the form 1/$k$ for some integer $k$, $\omega$ is 1 or
$-1$, which is impossible by Lemma \ref{zeroone}. Hence $\theta \neq
-4$.

To complete the proof of the lemma, it remains to show that
$\widetilde{W}(\ast)$ is hyperbolic.

$\widetilde{W}(\ast)$ is irreducible because $\widetilde{W}(1/0)$
and $\widetilde{W}(-1)$ are prime and are not homeomorphic. The way
of proving that $\widetilde{W}(\ast)$ is $\partial$-irreducible and
not a Seifert fibered space is exactly the same as for the case
$\varphi = -3$.

Suppose $\widetilde{W}(\ast)$ is toroidal. Let $T$ be an essential
torus in $\widetilde{W}(\ast)$. Since $\widetilde{W}(1/0)$ is a lens
space, $T$ must be separating. Let $\widetilde{W}(\ast)$ = $A \cup_T
B(\ast)$ where $B$ contains the torus boundary of
$\widetilde{W}(\ast)$.

\begin{claim}
$T$ is compressible in $\widetilde{W}(-1)$.
\end{claim}
\begin{proof}
Suppose $T$ is incompressible in $\widetilde{W}(-1)$. Let $T'$ be
the double branched cover of $S$ as shown in Figure 38.
$\widetilde{W}(-1)$ = $D^2(2, 2)$ $\cup_{T'}$ $D^2(2,
|1/\omega-1|)$. Then $T'$ is the unique incompressible torus up to
isotopy in $\widetilde{W}(-1)$ unless $\omega \in \{ -1, 0, 1, 1/0,$
$1/2, 1/3\}$. If $\omega = 1/3$, then $\widetilde{W}(1)$ is
reducible from Figure 38. Thus by Lemma \ref{rr}, $\Delta(\theta, 1)
\leq 1$ i.e. $\theta = 2$ since $\theta \neq 0, 1$, in which case we
are done. Therefore we may assume that $T'$ is the unique
incompressible torus up to isotopy in $\widetilde{W}(-1)$. Thus, $T
= T'$. In other words, $\widetilde{W}(-1)$ = $D^2(2, 2)$ $\cup_T$
$D^2(2, |1/\omega-1|)$. Hence $A$ is either $D^2(2, 2)$ or $D^2(2,
|1/\omega - 1)$. Observe from Figure 38 that $\widetilde{W}(0)$ =
$S^2(3, 3, |1/\omega-2|)$, $\widetilde{W}(1)$ = $S^2(4, 2,
|1/\omega-3|)$ and $\widetilde{W}(2)$ = $S^2(3, 2, |1/\omega-5|)$.
Hence $A$ is $D^2(2, 3)$ and $\omega$ = 1/4. $A$ = $D^2(2, 3)$ is
the complement of the trefoil knot i.e. the (2, 3) torus knot
$T_{2,3}$. Now apply the same argument as in the case $\varphi = -3$
using $\widetilde{W}(1)$ = $L(10, -3)$.
\end{proof}

$T$ is compressible in $\widetilde{W}(1/0)$, $\widetilde{W}(0)$,
$\widetilde{W}(1)$, $\widetilde{W}(-1)$ and $\widetilde{W}(2)$. Then
the rest of the proof is exactly the same as for the case $\varphi =
-3$.
\end{proof}

\begin{lemma} \label{-1/3}
If $\omega =$ $-1/3$, then $\theta =$ $2$.
\end{lemma}
\begin{proof}
Assume that $\omega = -1/3$. Let $W(\ast) = Q(\theta, \ast, -1/3,
-1/3)$. Then $\widetilde{W}(\varphi)$ is reducible. Consider the
rational tangle fillings shown in Figure 40 with the $-1/3$-rational
tangle substituted for $\omega$. Then $\widetilde{W}(-1)$ is
toroidal unless $\theta = 1/0, 0, 1, 2$ (by Lemma \ref{zeroone},
$\theta = 2$, in which case we are done.). Hence we may assume that
$\widetilde{W}(-1)$ is toroidal.

Assume $\widetilde{W}(\ast)$ is hyperbolic. $\widetilde{W}(\varphi)$
is reducible and $\widetilde{W}(-1)$ is toroidal. Thus
$\Delta(\varphi, -1) \leq 3$. By Lemma \ref{zeroone} $\varphi = -4,
-3, -2$ or 2. If $\varphi = -3$, then $\theta = 2$ by Lemma \ref{-3}.
If $\varphi = 2$, we are done by using the symmetry by rotation in
the horizontal axis on $Q$. Also, $\varphi \neq -2$
by Lemma \ref{-2}.

Suppose $\varphi = -4$. Let $W'(\ast) = Q(\ast, -4, -1/3, -1/3)$.
Then $\widetilde{W}'(\theta)$ is reducible. Similarly, we can
easily see from Figure 38 with $W$ replaced by $W'$ that $\widetilde{W}'(2)$ is the lens
space $L(5, 1)$ and $\widetilde{W}'(0)$ is $S^2(3, 4, 5)$. These
are non-homeomorphic prime manifolds. Thus $\widetilde{W}'(\ast)$ is
irreducible. Since $\widetilde{W}'(2)$ is a lens space and
$\widetilde{W}'(\theta)$ is reducible, by Lemma \ref{rr},
$\Delta(\theta, 2) \leq 1$. Thus $\theta$ is 1, 2 or 3. However 1 is
ruled out by Lemma \ref{zeroone}. If $\theta = 2$, then we are done.
The case $\theta = 3$ belongs to the case $\varphi=3$ by the symmetry
on $Q$. Thus $\theta = 2$ by Lemma \ref{3}.

To complete the proof of the lemma, it remains to show that
$\widetilde{W}(\ast)$ is hyperbolic.

$\widetilde{W}(\ast)$ is irreducible because $\widetilde{W}(1/0)$
and $\widetilde{W}(-1)$ are prime and not homeomorphic. In the proof
of Lemma \ref{-3}, in order to prove that $\widetilde{W}(\ast)$ is
not a Seifert fibered space, we used the fact that
$\widetilde{W}(-1)$ is irreducible and contains a separating
incompressible torus. In this case, $\widetilde{W}(-1)$ is also
irreducible and contains a separating incompressible torus unless
$\theta = -1, 0, 1, 2$, in which cases we are done by Lemma
\ref{zeroone}. Now we can apply the same argument as in the case
$\varphi = -3$ to show that $\widetilde{W}(\ast)$ is
$\partial$-irreducible and not a Seifert fibered space.

\begin{figure}[t]
\includegraphics{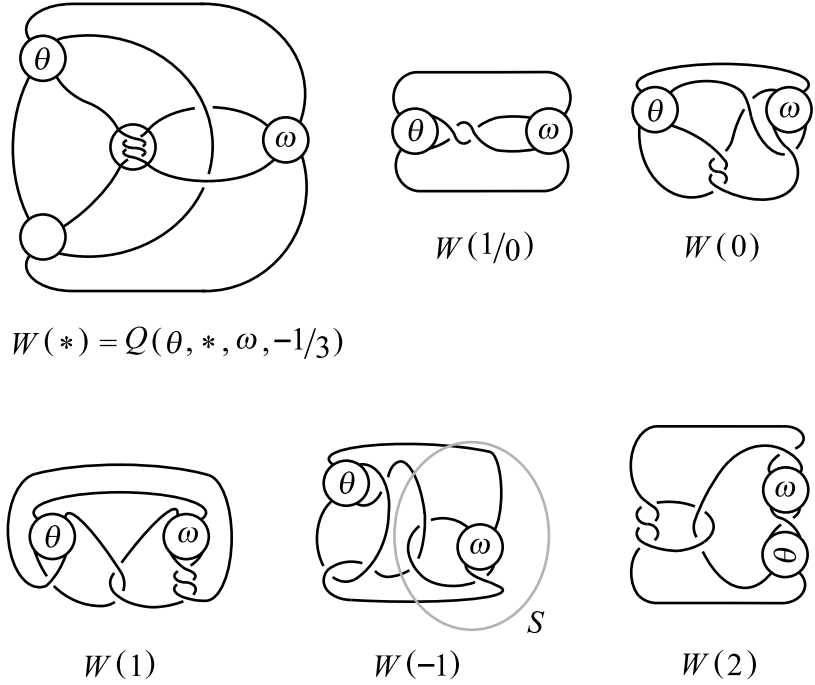}\caption{$W(\ast) = Q(\theta, \ast, \omega,
-1/3)$ and the rational tangle fillings.}\label{f:jj}
\end{figure}

Suppose $\widetilde{W}(\ast)$ is toroidal. Let $T$ be an essential
torus in $\widetilde{W}(\ast)$. Since $\widetilde{W}(1/0)$ is a lens
space, $T$ must be separating. Let $\widetilde{W}(\ast)$ = $A \cup_T
B(\ast)$ where $B$ contains the torus boundary of
$\widetilde{W}(\ast)$.

\begin{claim}
$T$ is compressible in $\widetilde{W}(-1)$.
\end{claim}
\begin{proof}
Suppose $T$ is incompressible in $\widetilde{W}(-1)$. Let $T'$ be
the double branched cover of $S$ as shown in Figure 40.
$\widetilde{W}(-1)$ = $D^2(2, 4)$ $\cup_{T'}$ $D^2(2,
|\theta-1|)$. Then $T'$ is the unique incompressible torus up to
isotopy in $\widetilde{W}(-1)$ unless $\theta = -1, 0, 1, 1/0, 2$,
in which cases we are done by Lemma \ref{zeroone}. Therefore $T' =
T$. That is, $\widetilde{W}(-1)$ = $D^2(2, 4)$ $\cup_T$ $D^2(2,
|\theta - 1|)$. Thus $A$ is either $D^2(2, 4)$ or $D^2(2,
|\theta - 1|)$. Observe from Figure 40 that $\widetilde{W}(0)$ =
$S^2(3, 5, |\theta|)$, $\widetilde{W}(1)$ = $S^2(2, 6, |\theta +
1|)$ and $\widetilde{W}(2)$ = $S^2(3, 2, |\theta + 5|)$. Hence $A$
is $D^2(2, 3)$ and $\theta$ = $-2$. $A$ = $D^2(2, 3)$ is the
complement of the trefoil knot i.e. the $(2,3)$ torus knot $T_{2,3}$.
Now apply the same argument as in the case $\varphi = -3$ using
$\widetilde{W}(1)$ = $L(16, -5)$.
\end{proof}

$T$ is compressible in $\widetilde{W}(1/0)$, $\widetilde{W}(0)$,
$\widetilde{W}(1)$, $\widetilde{W}(-1)$ and $\widetilde{W}(2)$. Then
the rest of the proof is exactly the same as for the case $\varphi = -3$.
\end{proof}

\begin{lemma} \label{1/3}
$\omega \neq$ $1/3$.
\end{lemma}
\begin{proof}
Suppose that $\omega= 1/3$. Let $W(\ast) = Q(\theta, \ast, 1/3,
-1/3)$. Then $\widetilde{W}(\varphi)$ is reducible. Consider the
rational tangle fillings shown in Figure 40 with the 1/3-rational
tangle substituted for $\omega$. Note that $\widetilde{W}(0)$ is a
lens space and $\widetilde{W}(1)$ is either a reducible manifold or
a lens space. Therefore Lemma \ref{rr} says that $\Delta(\varphi, 0)
\leq 1$ and $\Delta(\varphi, 1) \leq 1$ provided that $\widetilde{W}(\ast)$ is
irreducible. Hence $\varphi$ is 0 or 1 which is impossible by Lemma
\ref{zeroone}. However $\widetilde{W}(0)$ and $\widetilde{W}(-1)$
are prime manifolds and are not homeomorphic. This implies that
$\widetilde{W}(\ast)$ is irreducible.
\end{proof}

\begin{proof}[Proof of Proposition $\ref{phi}$]
Theorem \ref{23} says that the toroidal manifold $M(r_2)$ is either
$D^2(2, q)$ $\cup_T$ $D^2(2, s)$, $D^2(2, 3)$ $\cup_T$ $D^2(r, s)$
or $D^2(2, q)$ $\cup_T$ $D^2(3, s)$. On the other hand, $M(r_2)$ is
the double branched cover $\widetilde{Q}$ $( \theta, \varphi,
\omega, 1/0)$ over the link $Q( \theta, \varphi, \omega, 1/0)$ as
shown in Figure 41. Let $S$ be the sphere shown in Figure 41. The
double branched cover of $S$ is an incompressible torus in $M(r_2)$.
Since there is the unique incompressible torus in $M(r_2)$ up to
isotopy, $T$ is the double branched cover of $S$. Therefore, by
Lemma \ref{integers} $M(r_2)$ = $D^2(|\theta|, |\varphi|)$ $\cup_T$
$D^2(2, |1/\omega|)$.

\begin{figure}[tbh]
\includegraphics{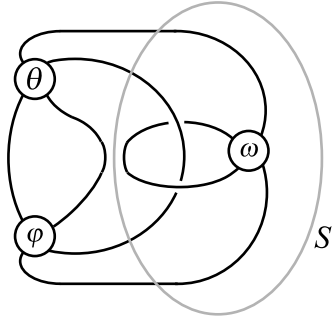}\caption{$Q( \theta, \varphi, \omega, 1/0)$.}\label{f:jj}
\end{figure}

(1) Assume $M(r_2)$ = $D^2(2, q)$ $\cup_T$ $D^2(2, s)$. Then
$\theta$ must be 2 or $-2$. Lemma \ref{-2} implies that $\theta$ is
2.

(2) Assume $M(r_2)$ = $D^2(2, q)$ $\cup_T$ $D^2(3, s)$. Then there
are two cases: either $\{|\theta|, |\varphi|\}$ = $\{3, s\}$ or
$\{2, |1/\omega|\}$ = $\{3, s\}$.

If $\{|\theta|, |\varphi|\}$ = $\{3, s\}$, then up to the symmetry,
we may assume that $|\varphi|$ = 3. Therefore by Lemmas \ref{-3},
\ref{3}, $\theta$ is 2.

If $\{2, |1/\omega|\}$ = $\{3, s\}$, then $|1/\omega|$ = 3. Lemmas
\ref{-1/3}, \ref{1/3} imply that $\theta$ is 2.

(3)Assume $M(r_2)$ = $D^2(2, 3)$ $\cup_T$ $D^2(r, s)$. There are two
cases: $\{|\theta|, |\varphi|\}$ = $\{2, 3\}$ or $\{2,
|1/\omega|\}$ = $\{2, 3\}$. In either case, Lemmas \ref{-3}, \ref{3},
\ref{-1/3} and \ref{1/3} show that $\theta$ is 2.
\end{proof}

\begin{proof}[Proof of Theorem $\ref{main2}$]
Proposition \ref{phi} says that $\theta$ = 2. Then $\widetilde{Q}(2,
\varphi, \omega,$ $-1/3)$ is $S^2(2, 3, |\varphi +2 - 1/\omega|)$. See $W(2)$ in Figure 38.
However since $M(r_1) = \widetilde{Q}(2, \varphi, \omega, -1/3)$ is
reducible, $\varphi +2 - 1/\omega$ must be 0. Thus if we let $\omega
= 1/p$ where $p$ is an integer, then $\varphi=p-2$. Hence $M(r_2) =
\widetilde{Q}(2, p-2, 1/p, -1/3)$. However Lemma \ref{zeroone}
implies that $p \neq -1, 0, 1, 2, 3$. If $p$ is either $-2$ or 4,
then $M(r_2) = \widetilde{Q}(2, p-2, 1/p, -1/3)$ contains a Klein
bottle, which is a contradiction. Therefore $p \leq -3$ or $p \geq
5$. Observe that the tangle $Q(2, p-2, 1/p, -1/3)$ is isotopic to the tangle
$Q(2, -p, 1/(-p+2), -1/3)$. Therefore, we may assume that $p \leq -3$.

To complete the proof of the theorem, we need only to show that
$\widetilde{Q}(2, p-2, 1/p, \ast)$ is hyperbolic. However, the proof
of this is exactly the same as the proof of Theorem 4.2 in \cite{EW}.

\end{proof}

\textbf{Acknowledgements}

The author would like to thank Sangyop Lee for his useful conversations.
Also, the author is very grateful to his advisor, Cameron McA Gordon, for his excellent suggestions and
encouragements through this work. We should also mention that since the author finished
writing this paper, he has learned that the main result of this paper has also recently been proved
by Sangyop Lee independently.


\begin{thebibliography}{CGLS}

\bibitem[1]{CGLS} M. Culler, C. McA. Gordon, J. Luecke, and P.B. Shalen,
            {\em Dehn surgery on knots}, Ann. of Math.
            \textbf{125} (1987), 237--300.
\bibitem[2]{EW} M. Eudave-Mu\~{n}oz and Y.Q. Wu,
            {\em Nonhyperbolic Dehn fillings on hyperbolic
            $3$-manifolds}, Pacific J.Math. \textbf{190} (1999),
            261--275.
\bibitem[3]{G1} C. McA. Gordon,
            {\em Boundary slopes on punctured tori in
            $3$-manifolds}, Trans. Amer. Math. Soc. \textbf{350}
            (1998), 1713--1790.
\bibitem[4]{G2} C. McA. Gordon,
            {\em Dehn surgery and satellite knots}, Trans. Amer. Math. Soc. \textbf{275}
            (1983), 687--708.
\bibitem[5]{GL} C. McA. Gordon and R.A. Litherland,
            {\em Incompressible planar surfaces in 3-manifolds},
            Topology Appl. \textbf{18} (1984), 121--144.
\bibitem[6]{GLu1} C. McA. Gordon and J. Luecke,
            {\em Knots are determined by their complements},
            J. Amer. Math. Soc. \textbf{2} (1989), 371--415.
\bibitem[7]{GLu2} C. McA. Gordon and J. Luecke,
            {\em Dehn surgeries on knots creating essential tori,
            I}, Communications in Analysis and Geometry
            \textbf{3} (1995), 597--644.
\bibitem[8]{GLu3} C. McA. Gordon and J. Luecke,
            {\em Reducible manifolds and Dehn surgery},
            Topology \textbf{35} (1996), 385--409.
\bibitem[9]{GLu4} C. McA. Gordon and J. Luecke,
            {\em Toroidal and boundary-reducing Dehn fillings},
            Topology Appl. \textbf{93} (1999), 77--90.
\bibitem[10]{GL4} C. McA. Gordon and J. Luecke,
            {\em Non-integral toroidal Dehn surgeries},
            Communications in Analysis and Geometry
            \textbf{12} (2004), 417--485.
\bibitem[11]{J} W. Jaco,
            {\em Lectures on three-manifold topology},
            Regional Conf. Series in Math. \textbf{43}(1997).
\bibitem[12]{Lee} S. Lee,
            {\em Reducing and Toroidal Dehn fillings on 3-manifold bounded by two tori},
            Mathematical Research Letters \textbf{12}(2005)
            10001--10014.
\bibitem[13]{Lee2} S. Lee,
            {\em Toroidal Dehn surgeries on Knots in $S^1 \times S^2$},
            J. Knot Theory Ramifications \textbf{14}(2005)
            657--664.
\bibitem[14]{LOT} S. Lee, S. Oh  and M. Teragaito,
            {\em Dehn fillings and small surfaces}, Preprint 2003.
            http://de.arxiv.org/abs/math.GT/0303156.
\bibitem[15]{O1} S. Oh,
            {\em Reducible and toroidal manifolds obtained by Dehn
            fillings}, Topology Appl. \textbf{75} (1997), 93--104.
\bibitem[16]{W1} Y.Q. Wu,
            {\em The reducibility of surgered $3$-manifolds},
            Topology Appl. \textbf{43} (1992), 213--218.
\bibitem[17]{W3} Y.Q. Wu,
            {\em Dehn fillings producing reducible manifolds and toroidal
            manifolds}, Topology \textbf{37} (1998), 95--108.
\end{thebibliography}
\end{document}